\documentclass[preprint,12pt]{elsarticle}

\usepackage{amsmath,amsfonts,amssymb,amsthm}
\usepackage{booktabs}
\usepackage{enumitem}
\usepackage{mathtools}
\usepackage{float}

\providecommand{\texorpdfstring}[2]{#1}

\newcommand{\R}{\mathbb{R}}
\newcommand{\rhoT}{\rho_\theta}

\newcommand{\hatG}{\widehat{G}}
\newcommand{\FQ}{F_Q}
\newcommand{\FV}{F_V}
\newcommand{\FR}{F_R}
\newcommand{\norm}[1]{\|#1\|}

\newcommand{\Ghat}{\widehat{G}}
\newcommand{\refd}{\mu}

\newtheorem{remark}{Remark}

\makeatletter
\def\ps@pprintTitle{%
  \let\@oddhead\@empty
  \let\@evenhead\@empty
  \let\@oddfoot\@empty
  \let\@evenfoot\@oddfoot}
\makeatother

\begin{document}

\begin{frontmatter}

\title{Computing the Gross–Pitaevskii Ground State via Wasserstein Gradient Flow in Diffeomorphism Space}

\author[purdue]{Xiangxiong Zhang\corref{cor1}}
\ead{zhan1966@purdue.edu}
\cortext[cor1]{Corresponding author}

\author[gatech]{Haomin Zhou}
\ead{hmzhou@math.gatech.edu}

\address[purdue]{Department of Mathematics, Purdue University,
150 N.\ University Street, West Lafayette, IN 47907, USA}
\address[gatech]{School of Mathematics, Georgia Institute of Technology,
Atlanta, GA 30332, USA}

\begin{abstract}
We compute the ground state $u$ of the Gross--Pitaevskii equation (GPE) via
Wasserstein gradient descent in diffeomorphism space.  We represent the density $\rho=u^2$ as the push-forward
of a fixed reference measure through a parameterized transport map
$T_\theta$, realized by a boundary-preserving Neural ODE.  The
Wasserstein gradient flow on probability densities then lifts to
natural gradient descent in the finite-dimensional parameter space, with
metric tensor given by the pullback of the Wasserstein metric.  The
method is entirely mesh-free and preserves the unit-mass constraint
without normalization.  We present numerical experiments in dimensions $d=1$ to
$4$, with separable and non-separable potentials, using a product
transport map in the separable case and a full map with coupled components
in the non-separable one.  The tests suggest that the output of the
parameterized Wasserstein gradient flow (PWGF) can be used as an effective
warm start for conventional solvers such as the $H^1$ Sobolev gradient
flow.
\end{abstract}


\end{frontmatter}

\section{Introduction}
\label{sec:intro}

A standard mathematical model of the
equilibrium states in Bose–Einstein condensation (BEC) is to consider
the ground state of the Gross--Pitaevskii equation (GPE) on a bounded
domain $\Omega\subset\R^d$ with homogeneous Dirichlet conditions, which is the
minimizer of the energy
\begin{equation}\label{eq:gpevp-energy-intro}
  E(u) = \frac{1}{2}\int_\Omega |\nabla u|^2 + V|u|^2
  + \frac{\beta}{2}|u|^4\,dx
\end{equation}
over the unit sphere
$\mathbb{S} = \{u\in H^1_0(\Omega): \norm{u}_{L^2}=1\}$,
with an external trapping potential  $V\ge 0$.  Equivalently, the ground state satisfies the
nonlinear eigenvalue problem
\begin{equation}\label{eq:gpevp-intro}
  -\Delta u + V u + \beta |u|^2 u = \lambda u, \qquad
  u\big|_{\partial\Omega}=0, \qquad \norm{u}_{L^2}=1,
\end{equation}
where the eigenvalue $\lambda$ is the Lagrange multiplier enforcing the constraint.  Existence, uniqueness, and
strict positivity of the ground state for  the defocusing regime $\beta\ge 0$ are
classical~\cite{lieb2000,lieb2005,zhang2024sobolev}.  

The recent review by
Henning and Jarlebring~\cite{henning2025review} organizes the extensive
literature on~\eqref{eq:gpevp-intro} into two complementary paradigms.
The \emph{algebraic paradigm}~\cite[Section~4]{henning2025review}
discretizes first, producing a nonlinear eigenvalue problem
$A(\mathbf{v})\mathbf{v}=\lambda\mathbf{v}$ in $\R^n$, then applies nonlinear eigenvalue algorithms  such as Self Consistent Field.
The \emph{variational paradigm}~\cite[Section~5]{henning2025review}
works with the continuous energy, equipping $\mathbb{S}$ with a Sobolev
inner product, such as $L^2$, $H^1$, or the solution-adapted $a_u$, to derive
gradient flows with progressively stronger convergence
guarantees~\cite[Theorems~5.15, 5.25]{henning2025review}. The gradient flow methods can also be regarded as Riemannian
optimization methods on unit sphere~\cite{danaila2017}. Convergence rate is often related to the spectral gap 
of the differential operator~\cite{henning2023}.    The two
paradigms are connected, e.g., the $J$-type Newton method coincides with
Rayleigh quotient iteration~\cite[Section~5.6.2]{henning2025review}.
All methods in both paradigms rely on a spatial discretization of
$\Omega$, with representation cost $O(n^d)$ in the number of mesh points
and per-step cost dominated by an elliptic solve (for Sobolev flows) or  
solving linear systems (for algebraic methods).   

In this paper we compute GPE ground state by Wasserstein gradient descent
in diffeomorphism space, extending the parameterized Wasserstein
Hamiltonian flow (PWHF) framework of Liu et
al.~\cite{liu2023siam,liu2024pwhf} from the time-dependent Schr\"odinger
equation to the nonlinear eigenvalue problem~\eqref{eq:gpevp-intro}.
Since the ground state satisfies $u>0$, the Madelung transform
$\rho=u^2$ rewrites the GP energy as a functional on probability
densities $\rho\in\mathcal{P}(\Omega)$:
\begin{equation}\label{eq:energy-rho-intro}
  E(\rho) = \underbrace{\tfrac{1}{8}\!\int_\Omega |\nabla\log\rho|^2
  \rho\,dx}_{\FQ(\rho)}
  + \underbrace{\tfrac{1}{2}\!\int_\Omega V\rho\,dx}_{\FV(\rho)}
  + \underbrace{\tfrac{\beta}{4}\!\int_\Omega \rho^2\,dx}_{\FR(\rho)},
\end{equation}
where $\FQ$ is the Fisher information (\S\ref{sec:math:madelung}).
The constraint $\norm{u}_{L^2}=1$ becomes $\int\rho\,dx=1$, which is
automatic for probability densities. The gradient descent with respect to
the $2$-Wasserstein metric gives
\begin{equation}\label{eq:wgf-intro}
  \partial_\tau\rho
  = \nabla\cdot\!\bigl(\rho\,\nabla\tfrac{\delta E}{\delta\rho}\bigr),
\end{equation}
a continuity equation that preserves mass by construction, unlike the
classical methods which usually require explicit renormalization.

We represent the density as the push-forward
$\rho_\theta = (T_\theta)_\sharp\refd$ of a fixed reference density
$\refd$ through a parameterized diffeomorphism
$T_\theta\colon\Omega\to\Omega$, realized by a boundary-preserving
Neural ODE (\S\ref{sec:algorithm}).
One draws $N$ sample points
$z_i\sim\refd$ once; at each iteration the particle positions
$x_i=T_\theta(z_i)$ and density values
$\rho_\theta(x_i)=\refd(z_i)/|\det\nabla T_\theta(z_i)|$ are obtained
by evaluating the map (\S\ref{sec:math:pushforward}).  The method
requires neither a spatial mesh nor a basis expansion, and preserves
the unit-mass constraint without normalization.

The Wasserstein metric on $\mathcal{P}(\Omega)$ pulls back to the
parameter space $\Theta\subset\R^M$ via the metric tensor
(\S\ref{sec:math:param})
\[
  \hatG(\theta) = \tfrac{1}{N}\textstyle\sum_{i=1}^N
  (\partial_\theta T_\theta(z_i))^\top(\partial_\theta T_\theta(z_i))
  \in\mathbb{R}^{M\times M},
\]
and the Wasserstein gradient descent becomes natural gradient
descent~\cite{amari1998}:
\begin{equation}\label{eq:pwgf-update-intro}
  \theta^{k+1} = \theta^k
  - \alpha\,\hatG(\theta^k)^+\nabla_\theta E(\theta^k),
\end{equation}
where $(\cdot)^+$ denotes a regularized pseudoinverse.
The parameters $\theta$ evolve by a geometric flow on the Wasserstein
manifold, discretized in parameter space.  The Fisher information $\FQ$
is computed via an augmented ODE that tracks the log-determinant and its
spatial derivative analytically, avoiding second-order automatic
differentiation (\S\ref{sec:alg:ode}).  Eigenvalues are
recovered from $\lambda = 2E + 2\FR$.

This work is inspired by 
the time-dependent PWHF framework of~\cite{liu2023siam,liu2024pwhf} in which
 $(\theta,p)$ evolves as a Hamiltonian system, and $p$ encodes the
quantum phase $\Phi$.  For the ground state, there is no Hamiltonian, thus 
we name \eqref{eq:pwgf-update-intro} as \emph{PWGF} (parameterized Wasserstein gradient flow).

In this paper, we provide a computational realization of PWGF for GPE
ground states in $d=1,2,3$ for separable potentials, and in
$d=2,3,4$ for non-separable ones. In PWGF, there is an approximation error from using neural network parameters $\theta$ to parametrize the diffeomorphism $T_\theta$.
The PWGF \eqref{eq:pwgf-update-intro} also contains an ill conditioned matrix $\hatG$, which often causes difficulty for high precision computation.
We benchmark against exact solutions, the
finite difference ground state references
of~\cite{zhang2024sobolev}, or spectral element reference solutions computed in this paper, and propose a
hybrid strategy in which PWGF output initializes the $H^1$ Sobolev
gradient flow, reducing the initial energy gap compared to
standard initializations.  

For simplicity, we only consider the defocusing case $\beta\ge 0$ on bounded
domains with Dirichlet conditions.  Convergence analysis of the
parameterized Wasserstein gradient flow remains open, and the numerical
results presented here are empirical.
Section~\ref{sec:math} derives the mathematical framework,
Section~\ref{sec:algorithm} describes the algorithm, and
Section~\ref{sec:numerics} presents numerical experiments
and the warm-start comparison, and concluding remarks are given in Section~\ref{sec:conclusion}.

\section{Mathematical Framework}
\label{sec:math}

This section derives the parameterized Wasserstein gradient flow in
general dimension~$d$.  We reformulate the GP energy in density form
(\S\ref{sec:math:madelung}), compute first
variations (\S\ref{sec:math:variations}), derive the Wasserstein gradient
flow (\S\ref{sec:math:wgf}), lift it to diffeomorphism space
(\S\ref{sec:math:pushforward}), and reduce to finite dimensions by Neural
ODE parameterization (\S\ref{sec:math:param}).

\subsection{Madelung transform: from eigenfunction to density}
\label{sec:math:madelung}

The ground state $u$ is strictly positive in $\Omega$ for the following reasons.  The energy satisfies $E(u)\ge E(|u|)$, so the minimizer can be taken
nonnegative.  A nonnegative minimizer satisfies~\eqref{eq:gpevp-intro}
with the self-adjoint operator $A_{u^*}=-\Delta+V+\beta(u^*)^2$, whose
Green's function is strictly positive by the strong maximum principle.
The Krein--Rutman theorem then implies $u^*>0$ on $\Omega$;
see~\cite{lieb2000,lieb2005,zhang2024sobolev} and, at the discrete level,
\cite[Theorem~3.4]{zhang2024sobolev}.

Because $u>0$, we may write
\begin{equation}\label{eq:madelung}
  u = \sqrt{\rho}, \qquad \rho = u^2 > 0.
\end{equation}
The kinetic energy $\tfrac{1}{2}\int|\nabla u|^2\,dx$ transforms into
the Fisher information of $\rho$.

\paragraph{Kinetic energy as Fisher information}
By writing $u=\sqrt\rho$ and applying the chain rule, we get
\begin{equation}\label{eq:kinetic_fisher}
  \frac12\int_\Omega |\nabla u|^2\,dx
  = \frac12\int_\Omega |\nabla\sqrt\rho|^2\,dx
  = \frac{1}{8}\int_\Omega |\nabla\log\rho|^2\,\rho\,dx,
\end{equation}
where $F_Q(\rho) := \frac{1}{8}\int_\Omega |\nabla\log\rho|^2\,\rho\,dx$ is the Fisher
information.
 
By substituting~\eqref{eq:kinetic_fisher}  and using
$|u|^2 = \rho$, $|u|^4 = \rho^2$, we get \eqref{eq:energy-rho-intro}, where
$\FQ$ is the Fisher information, $\FV$ the potential energy, and
$\FR$ the interaction energy.
The density formulation~\eqref{eq:energy-rho-intro} also appears in Bao and
Ruan~\cite{bao2019}, who observe that $|\nabla\rho|^2/\rho$ in $\FQ$
becomes singular as $\rho\to 0^+$, a difficulty for any grid-based method
that evaluates $\rho$ at fixed spatial nodes near $\partial\Omega$.

\begin{remark}[Regularization in grid-based density methods]
\label{rem:eps-reg}
Bao and Ruan~\cite{bao2019} address the singularity in $\FQ$ by
$\varepsilon$-regularization: replace $\rho$ by $\rho+\varepsilon$ in the
denominator, giving the regularized energy
\[
  E^\varepsilon(\rho)
  = \int_\Omega \frac{|\nabla\rho|^2}{8(\rho+\varepsilon)}
  + \frac{1}{2}V\rho + \frac{\beta}{4}\rho^2\,dx.
\]
They prove that $E^\varepsilon$ is convex for $\beta\ge 0$, that the
minimizer $\rho_g^\varepsilon\to\rho_g$ as
$\varepsilon\to 0^+$. The
convergence rate $\|\rho_g - \rho_g^\varepsilon\|_{L^2} = O(\sqrt\varepsilon)$ was expected from their analysis, though
the asymptotic error $O(\varepsilon)$ can be observed numerically in some setup.   
\end{remark}

\paragraph{Constraint transformation}
The $L^2$-normalization $\norm{u}_{L^2}=1$ becomes $\int_\Omega\rho\,dx=1$:
the density is automatically a probability measure.  The Dirichlet boundary
condition $u|_{\partial\Omega}=0$ becomes $\rho|_{\partial\Omega}=0$.

\subsection{First variations and eigenvalue recovery}
\label{sec:math:variations}

We compute the first variation of each term in~\eqref{eq:energy-rho-intro} with respect
to $\rho$, subject to $\int\rho\,dx=1$.

\paragraph{First variation of $\FQ$}
The first variation is
\begin{equation}\label{eq:varFQ}
  \frac{\delta\FQ}{\delta\rho} = -\frac{1}{2}\frac{\Delta\sqrt\rho}{\sqrt\rho}.
\end{equation}
To verify, write $\FQ(\rho) = \tfrac{1}{2}\int|\nabla\sqrt\rho|^2\,dx$ and
take the G\^ateaux derivative: for a perturbation $\sigma$ with
$\int\sigma\,dx=0$,
\[
  \left.\frac{d}{d\varepsilon}\right|_{\varepsilon=0}
  \FQ(\rho+\varepsilon\sigma)
  = \frac{1}{2}\int_\Omega \nabla\sqrt\rho \cdot
    \nabla\!\left(\frac{\sigma}{2\sqrt\rho}\right)\,dx
  = -\frac{1}{2}\int_\Omega \frac{\Delta\sqrt\rho}{\sqrt\rho}\cdot
    \frac{\sigma}{2}\,dx + \text{(b.t.)}
\]
after integrating by parts (boundary terms vanish since $\rho|_{\partial\Omega}=0$).
Reading off the coefficient of $\sigma/2$
gives~\eqref{eq:varFQ}.  The quantity
$-\tfrac{1}{2}\Delta\sqrt\rho/\sqrt\rho$ is the \emph{quantum potential}
(Bohm potential), connecting the Madelung and de~Broglie--Bohm pictures.

\paragraph{First variations of $\FV$ and $\FR$}
The remaining variations are immediate:
\begin{equation}\label{eq:varFVFR}
  \frac{\delta\FV}{\delta\rho} = \frac{1}{2}V, \qquad
  \frac{\delta\FR}{\delta\rho} = \frac{\beta}{2}\rho.
\end{equation}

\paragraph{Euler--Lagrange equation}
At a minimizer of $E(\rho)$ subject to $\int\rho\,dx=1$, the first-order
optimality condition with Lagrange multiplier $\lambda$ reads
\begin{equation}\label{eq:EL_rho}
  \frac{\delta E}{\delta\rho}
  = -\frac{1}{2}\frac{\Delta\sqrt\rho}{\sqrt\rho}
    + \frac{1}{2}V + \frac{\beta}{2}\rho
  = \frac{\lambda}{2}.
\end{equation}
Multiplying both sides by $2\sqrt\rho$ and writing $u=\sqrt\rho$ recovers
the GPE~\eqref{eq:gpevp-intro}.

\paragraph{Eigenvalue formula}
Multiplying the GPE~\eqref{eq:gpevp-intro} by $u$ and integrating by parts,
\begin{equation}
  \lambda = \int_\Omega\bigl(|\nabla u|^2 + V u^2 + \beta u^4\bigr)\,dx
          = 2\FQ + 2\FV + 4\FR.
  \label{eq:lam_raw}
\end{equation}
Since $E = \FQ + \FV + \FR$, we obtain
\begin{equation}
  \lambda = 2E + 2\FR.
  \label{eq:lam_est}
\end{equation}
This is the eigenvalue estimator used in all our numerical experiments.

\subsection{The Wasserstein gradient flow}
\label{sec:math:wgf}

Having reformulated the ground state as the minimizer of $E(\rho)$ over
probability densities, we now choose a \emph{metric} on the space of densities
to define the gradient descent.  Different metrics for $u$ give different gradient flows, e.g.,
the $L^2$ metric leads to the normalized gradient flow of Bao and
Du~\cite{bao2004} and a Sobolev metric leads to the Sobolev gradient flow~\cite{henning2020}.  Here we use the 2-Wasserstein metric
from optimal transport for $\rho$.

\paragraph{The 2-Wasserstein distance}
For two probability densities $\rho_0,\rho_1\in\mathcal{P}(\Omega)$, the
2-Wasserstein distance is
\begin{equation}\label{eq:W2}
  W_2(\rho_0,\rho_1)^2
  = \inf_\pi \int_{\Omega\times\Omega} |x-y|^2\,d\pi(x,y),
\end{equation}
where the infimum is over all couplings $\pi$ with marginals $\rho_0$ and
$\rho_1$~\cite{villani2003,ambrosio2008}.  The Wasserstein distance metrizes
$\mathcal{P}(\Omega)$ with a geometry tied to the \emph{cost of transporting
mass}, rather than to pointwise comparison of densities.

\paragraph{Wasserstein metric tensor (Otto calculus)}
Following Otto~\cite{otto2001}, the tangent space at $\rho\in\mathcal{P}(\Omega)$
consists of perturbations $\sigma$ satisfying $\int\sigma\,dx=0$.  Each such
$\sigma$ is associated with a velocity potential $\Phi$ via the continuity
equation $\sigma = -\nabla\cdot(\rho\nabla\Phi)$.  The Wasserstein Riemannian
metric is
\begin{equation}\label{eq:gW}
  g^W_\rho(\sigma_1,\sigma_2)
  = \int_\Omega \nabla\Phi_1\cdot\nabla\Phi_2\;\rho\,dx,
  \qquad \text{where } -\nabla\cdot(\rho\nabla\Phi_i)=\sigma_i.
\end{equation}
This is an \emph{inner product on tangent directions} that weights
perturbations by how much mass must be moved.

\paragraph{Wasserstein gradient}
The Wasserstein gradient $\operatorname{grad}_W E(\rho)$ is the unique tangent
vector satisfying $g^W_\rho(\operatorname{grad}_W E,\,\sigma) = dE(\rho)[\sigma]$
for all admissible $\sigma$.  Using~\eqref{eq:gW} and integrating by parts,
one finds
\begin{equation}\label{eq:gradW}
  \operatorname{grad}_W E(\rho) = -\nabla\cdot\!\left(\rho\,\nabla
  \frac{\delta E}{\delta\rho}\right).
\end{equation}

\paragraph{The gradient flow}
Steepest descent $\partial_\tau\rho = -\operatorname{grad}_W E(\rho)$ gives
\begin{equation}
  \partial_\tau\rho = \nabla\cdot\!\left(\rho\,\nabla
  \frac{\delta E}{\delta\rho}\right).
  \label{eq:wgf}
\end{equation}
This PDE preserves mass ($\int\rho\,dx=1$) automatically because the
right-hand side is a divergence.  No normalization step is needed.

The stationary point of~\eqref{eq:wgf} satisfies
$\nabla(\delta E/\delta\rho)=0$, i.e.,
$\delta E/\delta\rho = \text{const} = \lambda/2$, recovering the
Euler--Lagrange equation~\eqref{eq:EL_rho} and hence the GPE.

\paragraph{Relation to Sobolev gradient flows}
In the $H^1$-gradient flow~\cite{zhang2024sobolev}, the
gradient step is $\partial_\tau u = -(-\Delta+a)^{-1}(\delta E/\delta u)$,
requiring an elliptic solve per step and explicit normalization.  The
Wasserstein flow~\eqref{eq:wgf} operates on $\rho$ rather than $u$,
preserves mass by construction, and its per-step cost depends on the
parameterization (see \S\ref{sec:math:param} below) rather than on a mesh-based
elliptic solver.

\paragraph{No momentum for the ground state}
For the time-dependent Schr\"odinger equation, the Madelung system involves
both $\rho$ and a phase $\Phi$, leading to a Wasserstein \emph{Hamiltonian}
flow~\cite{liu2023siam,liu2024pwhf}.  For the ground state, $\Phi\equiv 0$
(no phase), and no momentum variable $p$ is needed.  This is one particular difference compared to the time-dependent PWHF framework.

\subsection{Lifting to diffeomorphism space}
\label{sec:math:pushforward}

Following~\cite{liu2023siam,liu2024pwhf}, we lift the gradient
flow~\eqref{eq:wgf} from $\mathcal{P}(\Omega)$ to the space of
diffeomorphisms
$\mathcal{O}=\{T:\Omega\to\Omega : T\text{ is a diffeomorphism}\}$ by
representing each density as the push-forward of a fixed reference density
through a transport map~$T$.  The Wasserstein metric pulls back to the
$L^2(\refd)$-metric on $\mathcal{O}$ (see~\eqref{eq:gO}), and the
gradient flow becomes a pointwise ODE~\eqref{eq:flow_T} whose trajectory
is fully determined by its initial condition.  All integrals over $\rho$
become sample averages over fixed reference particles
$z_i\sim\refd$~\cite{liu2024pwhf}.

\medskip
\paragraph{Push-forward map}\;
Let $\refd$ be a fixed \emph{reference density} on $\Omega$ (chosen to
satisfy $\refd|_{\partial\Omega}=0$, with the same boundary vanishing as
$\rho$).  An orientation-preserving diffeomorphism $T:\Omega\to\Omega$ pushes
$\refd$ forward to
\begin{equation}\label{eq:pushforward}
  \rho = T_\sharp\refd, \qquad
  \rho(T(z)) = \frac{\refd(z)}{|\det\nabla T(z)|}
  \quad\text{for } z\in\Omega.
\end{equation}
By the change of variables formula, $\int_\Omega\rho\,dx=\int_\Omega\refd\,dz=1$,
so the push-forward preserves the probability constraint automatically.  Expectations
under $\rho$ become expectations under $\refd$:
\[
  \int_\Omega \varphi(x)\,\rho(x)\,dx = \int_\Omega \varphi(T(z))\,\refd(z)\,dz
  \approx \frac{1}{N}\sum_{i=1}^N \varphi(T(z_i)),
\]
where $z_1,\ldots,z_N\sim\refd$ are fixed quadrature particles.

\medskip
\paragraph{Metric on diffeomorphism space}\;
The Wasserstein metric~\eqref{eq:gW} on $\mathcal{P}(\Omega)$ induces a
metric on the space $\mathcal{O} = \{T:\Omega\to\Omega \text{ diffeo.}\}$.
For two tangent vectors $\sigma_1,\sigma_2\in T_T\mathcal{O}$ (vector fields
on $\Omega$), the induced metric is~\cite{liu2023siam}
\begin{equation}\label{eq:gO}
  \mathfrak{g}_T(\sigma_1,\sigma_2)
  = \int_\Omega \sigma_1(z)^\top \sigma_2(z)\;\refd(z)\,dz.
\end{equation}
This is simply the $L^2(\refd)$-inner product on vector fields.
This pullback is exact \cite{liu2023siam}.

\medskip
\paragraph{Gradient flow on $\mathcal{O}$}\;
The Wasserstein gradient descent~\eqref{eq:wgf} pulls back to
\begin{equation}\label{eq:flow_T}
  \partial_\tau T(z) = -\nabla_x \frac{\delta E}{\delta\rho}
  \big(T_\sharp\refd,\,\cdot\,\big)\circ T(z).
\end{equation}
This is an ODE on the space of diffeomorphisms, driven by the gradient of
$\delta E/\delta\rho$ evaluated along the current map.

\subsection{Finite-dimensional reduction via Neural ODE}
\label{sec:math:param}
 
We restrict $T$ to a finite-dimensional family $\mathcal{O}_\Theta =
\{T_\theta:\theta\in\Theta\subset\R^M\}$, where $T_\theta$ is a Neural ODE
(details in Section~\ref{sec:alg:ode}).

\medskip
\paragraph{Pullback metric tensor}\;
A tangent vector $\dot\theta\in\R^M$ in parameter space maps to a tangent
vector $\sigma = (\partial_\theta T_\theta)\dot\theta$ in diffeomorphism
space.  Substituting into the metric~\eqref{eq:gO}:
\begin{align*}
  \mathfrak{g}(\dot\theta_1,\dot\theta_2)
  &= \int_\Omega \bigl[(\partial_\theta T_\theta)\dot\theta_1\bigr]^\top
     \bigl[(\partial_\theta T_\theta)\dot\theta_2\bigr]\;\refd\,dz \\
  &= \dot\theta_1^\top
     \underbrace{\left[\int_\Omega
       (\partial_\theta T_\theta(z))^\top (\partial_\theta T_\theta(z))
       \;\refd(z)\,dz
     \right]}_{G(\theta)}\,
     \dot\theta_2.
\end{align*}
Here $\partial_\theta T_\theta(z)\in\R^{d\times M}$ is the Jacobian of
$T_\theta(z)$ with respect to $\theta$, so
$G(\theta)\in\R^{M\times M}$ is the \emph{pullback metric tensor}.
Approximating the $\refd$-integral by Monte Carlo with the fixed particles
$z_1,\ldots,z_N\sim\refd$:
\begin{equation}
  \Ghat(\theta) = \frac{1}{N}\sum_{i=1}^{N}
    \bigl(\partial_\theta T_\theta(z_i)\bigr)^\top
    \bigl(\partial_\theta T_\theta(z_i)\bigr)
    \;\in\;\R^{M\times M}.
  \label{eq:G}
\end{equation}

\medskip
\paragraph{Natural gradient descent}\;
Restricting the gradient flow~\eqref{eq:flow_T} to $\mathcal{O}_\Theta$ and
using the pullback metric, the Wasserstein gradient descent in parameter space
becomes
\begin{equation}
  \theta^{k+1} = \theta^k - \alpha\,\Ghat(\theta^k)^+\,\nabla_\theta E(\theta^k),
  \label{eq:wgd}
\end{equation}
where $(\cdot)^+$ denotes a regularized pseudoinverse and $\alpha>0$ is the
step size.  This is \emph{natural gradient descent}~\cite{amari1998}: the
matrix $\Ghat(\theta)^+$ rotates and rescales the Euclidean gradient
$\nabla_\theta E$ to account for the non-Euclidean geometry of the
parameter-to-density map $\theta\mapsto\rhoT$.   

\begin{remark}[No regularization needed in PWGF]
\label{rem:no-reg}
Unlike grid-based density methods~\cite{bao2019}, the push-forward
parameterization avoids the singularity in $\FQ$ without any
$\varepsilon$-regularization.  In PWGF, the density $\rhoT$ is never
evaluated at arbitrary grid points; it is evaluated only at the pushed
particle locations $x_i = T_\theta(z_i)$, where
\[
  \rhoT(x_i) = \frac{\refd(z_i)}{|\det\nabla T_\theta(z_i)|}.
\]
This ratio is strictly positive because:
\emph{(i)}~the particles $z_i$ are sampled from the interior of
$\Omega$, where the reference density satisfies $\refd(z_i)>0$; and
\emph{(ii)}~$T_\theta$ is an orientation-preserving diffeomorphism, so
$\det\nabla T_\theta(z_i)>0$.
The score $\partial_x\log\rhoT$ is computed   via the
augmented ODE (see Section~\ref{sec:alg:ode}), not by numerical
differentiation of $\rho$, eliminating the $|\nabla\rho|^2/\rho$
singularity at the discrete level.
The only requirements are that $\refd$ vanish at $\partial\Omega$
and that sample points avoid the boundary; both hold when $z_i$ are
drawn from a Beta distribution on the open interval.  
\end{remark}

\subsection{Summary of PWGF}
\label{sec:math:summary}

We summarize the PWGF approach and compare it with two families of
gradient-flow methods for the GPE ground state.

\medskip
\paragraph{Conceptual steps of PWGF}\;
For the GPE ground state eigenvalue problem~\eqref{eq:gpevp-intro}:
\begin{enumerate}[label=(\roman*)]
  \item \emph{Madelung transform.}  Use $u>0$ to write $\rho=u^2$,
        converting the energy $E(u)$ into a functional
        $E(\rho) = \FQ(\rho)+\FV(\rho)+\FR(\rho)$ on probability
        densities (\S\ref{sec:math:madelung}).
  \item \emph{Wasserstein gradient.}  Choose the 2-Wasserstein metric on
        $\mathcal{P}(\Omega)$ and derive the gradient flow
        $\partial_\tau\rho = \nabla\cdot(\rho\,\nabla(\delta E/\delta\rho))$
        (\S\ref{sec:math:wgf}).
  \item \emph{Lift to diffeomorphism space.}  Represent $\rho=T_\sharp\refd$
        and pull back the gradient flow to an ODE on the space of maps:
        $\partial_\tau T = -\nabla_x(\delta E/\delta\rho)\circ T$
        (\S\ref{sec:math:pushforward}).
  \item \emph{Parameterize and reduce.}  Restrict $T$ to a Neural ODE
        family $T_\theta$, yielding natural gradient descent
        $\theta^{k+1} = \theta^k - \alpha\,\Ghat^+\nabla_\theta E$
        in $\R^M$ (\S\ref{sec:math:param}).
\end{enumerate}
Steps (iii)--(iv) are the distinctive features of PWGF: the lifting to
diffeomorphism space replaces the density-space PDE with a sampling-based
ODE, and the Neural ODE parameterization reduces the problem to
finite-dimensional natural gradient descent without spatial
discretization.

\medskip
\paragraph{Sobolev gradient flows.}\;
The Sobolev gradient flow framework~\cite{henning2020,zhang2024sobolev}
works directly with the eigenfunction $u$ rather than the
density $\rho$.  Given a Hilbert space $X$ with inner
product $(\cdot,\cdot)_X$, the $X$-Sobolev gradient $\nabla_X E(u)$ is
defined as the Riesz representative of $E'(u)$:
$(\nabla_X E(u),\,w)_X = \langle E'(u),\,w\rangle$ for all $w\in X$.
The gradient flow on the unit sphere $\mathbb{S}$ is then
$\partial_\tau u = -P_{u,X}\nabla_X E(u)$, where $P_{u,X}$ projects onto
the tangent space of $\mathbb{S}$ at $u$.  Different choices of $X$ give
different methods:
\begin{itemize}
  \item \emph{$L^2$ gradient with $(w,z)_X=(w,z)$} ($X=L^2$):
        $\nabla_{L^2}E(u) = -\Delta u + Vu + \beta u^3$, leading to the
        GFDN scheme(semi-implicit gradient flow)~\cite{bao2004};
  \item \emph{$H^1$ gradient} ($X=H^1_0$ with $(w,z)_X=(\nabla w,\nabla z)+(w,z)$):
        computing $(-\Delta+I)^{-1}$ per step,
        giving the $H^1$-flow~\cite{zhang2024sobolev};
  \item \emph{$a_0$ gradient} ($X=H^1_0$ with $(w,z)_X=(\nabla w,\nabla z)+(Vw,z)$):
        computing $(-\Delta+V)^{-1}$ per step, admits adaptive step
        sizes~\cite{henning2025review};
  \item \emph{$a_u$ gradient} ($X=H^1_0$ with
        $(w,z)_X=(\nabla w,\nabla z)+(Vw,z)+\beta(u^2 w,z)$):
        $\nabla_{a_u}E(u)=u$ trivially, reducing the flow to damped
        inverse iteration~\cite{henning2025review}.
\end{itemize}
All of these are mesh-based: they require a spatial discretization
(finite element or finite difference) and solve one or two elliptic PDEs per
iteration.  They enjoy rigorous convergence guarantees
such as global energy decrease and local exponential convergence with rate governed
by the spectral gap of the differential operator involved. See~\cite{henning2025review} and references therein.

By contrast, PWGF is mesh-free: it uses $N$ fixed particles and a Neural
ODE with $M$ parameters, and computes the natural gradient direction by
solving an $M\times M$ linear system rather than an elliptic PDE.  The
Wasserstein metric automatically preserves the mass constraint
$\int\rho\,dx=1$, so no normalization step is needed, whereas all Sobolev
flows require explicit projection back to $\mathbb{S}$.

\medskip
\paragraph{Direct Wasserstein discretizations}\;
The Wasserstein gradient flow~\eqref{eq:wgf} on $\mathcal{P}(\Omega)$ is
the mathematical starting point for PWGF, but discretizing it directly
(e.g., via the JKO scheme~\cite{jordan1998})
requires a spatial mesh and suffers from the curse of dimensionality.
The PWGF approach differs from such direct discretizations in two ways:
\begin{enumerate}[label=(\alph*)]
  \item \emph{Lifting.}  Rather than discretizing the density-space
        PDE~\eqref{eq:wgf}, PWGF lifts it to the space of
        diffeomorphisms~$\mathcal{O}$, where the flow becomes a pointwise
        ODE~\eqref{eq:flow_T}.  This ODE is self-contained: the
        trajectory of $T$ is fully determined by its initial condition,
        without solving the continuity equation separately.
  \item \emph{Parameterization.}  The Neural ODE $T_\theta$ reduces the
        infinite-dimensional ODE on $\mathcal{O}$ to a finite-dimensional
        system in~$\R^M$.  The pullback of the Wasserstein metric gives the
        metric tensor~$\Ghat$, and the flow becomes natural gradient
        descent~\eqref{eq:wgd}.
\end{enumerate}
The price is the difficulty of rigorous convergence guarantees. In PWGF, the energy is not
guaranteed to decrease monotonically since the finite dimensional parameterization introduces
approximation error and the linear solver is only an approximation.  A rigorous convergence proof of PWGF is widely open.

\medskip
\paragraph{Related work}\;
The parameterized Wasserstein gradient flow framework was introduced
in~\cite{jin2025pwgf} for classical Wasserstein gradient flows
(Fokker--Planck, porous medium, and aggregation equations), and its
Hamiltonian counterpart was developed in~\cite{liu2023siam} for the
time-dependent Schr\"odinger equation.
Neither paper addresses the GPE \emph{ground state} (nonlinear eigenvalue)
problem, which is the focus of the present work.

To the best of our knowledge, no prior work applies Wasserstein gradient
flow to the GPE eigenvalue problem.
The $L^2$ normalized gradient flow for the GPE ground state has been
analyzed extensively (see~\cite{bao2004} and the recent
work~\cite{chu2025gfgpe}); Sobolev gradient flows ($H^1$, $a_0$, $a_u$
metrics) are treated in~\cite{henning2020,henning2023,zhang2024sobolev,
henning2025review}.
All of these operate in function space and require spatial discretization.
In a different direction, Bao et al.~\cite{bao2025normdnn} minimize the
GP energy via a normalized deep neural network with standard training;
their method is not based on Wasserstein geometry or push-forward maps.
The push-forward parameterization $\rho_\theta = (T_\theta)_\sharp\refd$ is
closely related to normalizing flows~\cite{rezende2015,papamakarios2021};
the key difference is that PWGF evolves parameters by a geometric ODE rather
than by maximum-likelihood training.
Neklyudov et~al.~\cite{neklyudov2023wqmc} use a Wasserstein gradient
flow on the Born distribution space to accelerate quantum variational
Monte Carlo for fermionic systems, but work directly in the space of
distributions (without lifting to diffeomorphism space) and target the
many-body Schr\"odinger equation rather than the mean-field GPE.
Peterseim et~al.~\cite{peterseim2025nn} use neural networks to
accelerate conventional iterative solvers for nonlinear Schr\"odinger
eigenvalue problems, leveraging knowledge from prior simulations to
predict improved solution trajectories.

\section{Algorithm}
\label{sec:algorithm}

Implementing the natural gradient descent~\eqref{eq:wgd} requires three
ingredients: a parameterization $T_\theta$ that evaluates the density
$\rhoT$ and score $\nabla_x\log\rhoT$
(\S\ref{sec:alg:ode}); a discretization of $E(\theta)$ and its gradient
(\S\ref{sec:alg:energy}); and a solver for the $M\times M$ system
$\Ghat\xi=\nabla_\theta E$ (\S\ref{sec:alg:solver}).  The complete
algorithm is stated in \S\ref{sec:alg:main}; dimension-specific choices
are in \S\ref{sec:alg:dims}.

\subsection{Boundary-preserving Neural ODE}
\label{sec:alg:ode}

On $(-L,L)^d$ with Dirichlet conditions, the boundary constraint
$\rhoT|_{\partial\Omega}=0$ requires $T_\theta(\pm L)=\pm L$.  We
enforce this by defining the Neural ODE velocity field (per coordinate)
as
\begin{equation}
  f_\theta(w) = \bigl(1-w^2/L^2\bigr)\,g_\theta(w),
  \label{eq:bpode}
\end{equation}
where $g_\theta:\R\to\R$ is a tanh network $1\!\to\!H\!\to\!H\!\to\!1$
with $M_1 = H^2+4H+1$ parameters.
The factor $(1-w^2/L^2)$ vanishes at $w=\pm L$, so the ODE
$\dot w = f_\theta(w)$ has fixed points at $\pm L$: starting from any
$w(0)=z\in(-L,L)$, the trajectory remains in $(-L,L)$ and
$T_\theta(\pm L)=\pm L$ for all~$\theta$.

\paragraph{Spatial derivatives}
The energy computation requires the first and second derivatives of
$f_\theta$ with respect to the spatial variable~$w$:
\begin{align}
  f'_\theta(w)  &= -\frac{2w}{L^2}\,g_\theta(w)
                    + \bigl(1-w^2/L^2\bigr)\,g'_\theta(w),
  \label{eq:fp}\\
  f''_\theta(w) &= -\frac{2}{L^2}\,g_\theta(w)
                    - \frac{4w}{L^2}\,g'_\theta(w)
                    + \bigl(1-w^2/L^2\bigr)\,g''_\theta(w).
  \label{eq:fpp}
\end{align}
Here $g'_\theta$ and $g''_\theta$ are computed analytically via the chain
rule through the tanh layers.

\paragraph{Augmented ODE}
Following~\cite{liu2024pwhf}, we evolve an augmented state
$(w,\,\ell,\,J,\,\ell')$ with dynamics
\begin{align*}
  \dot w    &= f_\theta(w),\\
  \dot\ell  &= f'_\theta(w),\\
  J         &= \exp(\ell) \quad\text{(updated after each $\ell$ step)},\\
  \dot\ell' &= f''_\theta(w)\cdot J,
\end{align*}
where $\ell = \log|T_\theta'(z)|$ is the log-Jacobian determinant,
$J = |T_\theta'(z)| = e^\ell$ is the Jacobian determinant itself, and
$\ell' = \partial_z\log|T_\theta'(z)|$.  All four components are evolved
per coordinate for each of the $N$ particles, using forward Euler with
$N_{\rm ODE}$ substeps.

\paragraph{Score function}
 Through the ODE, the augmented state avoids second-order differentiation,
which causes severe numerical instabilities~\cite{liu2024pwhf}.
The score at $x_k = T_{\theta_k}(z_k)$ for coordinate~$k$ is
\begin{equation}
  \partial_{x_k}\log\rhoT
  = \frac{\partial_{z_k}\log\refd_k(z_k) - \ell'_k}{T_{\theta_k}'(z_k)},
  \label{eq:score}
\end{equation}
where $\partial_{z_k}\log\refd_k$ is the score of the reference density
(Section~\ref{sec:alg:dims}).

\subsection{Energy and its gradient}
\label{sec:alg:energy}

Given $N$ particles $z_i\sim\refd$ and the augmented ODE outputs at
$x_i=T_\theta(z_i)$, the energy terms~\eqref{eq:energy-rho-intro} become sample
averages via the change of variables
$\int h(x)\rho\,dx = \int h(T(z))\refd\,dz$:
\begin{align}
  \FQ(\theta) &\approx \frac{1}{8N}\sum_{i=1}^N
    \sum_{k=1}^{d} s_{k,i}^2,
  \quad s_{k,i} = \frac{\partial_{z_k}\log\refd_k(z_{k,i})-\ell'_{k,i}}
       {J_{k,i}},
  \label{eq:FQ_disc}\\
  \FV(\theta) &\approx \frac{1}{2N}\sum_{i=1}^N V(x_i),
  \label{eq:FV_disc}\\
  \FR(\theta) &\approx \frac{\beta}{4N}\sum_{i=1}^N \rhoT(x_i),
  \quad \rhoT(x_i) = e^{\log\refd(z_i)-\ell_i},
  \label{eq:FR_disc}
\end{align}
where for the product-map ansatz (Section~\ref{sec:alg:dims}) the total
log-Jacobian determinant factorizes additively as
$\ell_i = \sum_{k=1}^d \ell_{k,i}$, and similarly
$\log\refd(z_i)=\sum_k\log\refd_k(z_{k,i})$.

The eigenvalue is estimated from~\eqref{eq:lam_est}: $\lambda = 2E + 2\FR$.

\paragraph{Gradient computation}
The gradient $\nabla_\theta E$ is obtained by automatic differentiation
through the forward ODE: all $N$ particles are integrated in a single
batched pass, and the computational graph is differentiated in the
standard reverse mode.
The metric tensor $\Ghat$, on the other hand, requires the per-sample
Jacobians $\partial_\theta T_\theta(z_i) \in \R^{d\times M}$, which are
computed via a separate forward pass using functional differentiation
(reverse-mode Jacobian computation vectorized over the $N$ samples),
cf.~\eqref{eq:G}.

\subsection{Linear solver for the natural gradient}
\label{sec:alg:solver}

The natural gradient direction $\xi$ solves
$\Ghat(\theta)\,\xi = \nabla_\theta E(\theta)$, where $\Ghat$ is
symmetric positive semi-definite. However, $\Ghat$ is highly ill conditioned with an initial condition number around
$O(10^{10})$ and $\Ghat$ is often rank deficient in our experiments.

We use conjugate gradients (CG) with Tikhonov regularization
$(\Ghat + \epsilon I)\xi = \nabla_\theta E$, $\epsilon = 10^{-6}$,
run for up to $n_{\rm CG}$ iterations (100 for 1D, 200 for 2D, 300 for
3D).  The regularization $\epsilon=10^{-6}$ is chosen to ensure the CG
iteration converges in a controlled number of steps.  After computing the CG solution $\xi$, a backtracking safety check
evaluates $E(\theta-\alpha\xi)$ via a trial forward pass; if
$E_{\rm new}>E+E_{\rm tol}$, the step is rejected in favor of a small
normalized gradient descent step.  Throughout the iteration, the
parameters achieving the lowest energy (\emph{best parameters}) are
recorded for use in the final reconstruction.

All steps are followed by norm clipping:
$\xi \leftarrow \min(1,\,C/\norm{\xi})\,\xi$.
Table~\ref{tab:hyperparams} lists the solver parameters for each problem.

\subsection{The complete algorithm}
\label{sec:alg:main}

The complete PWGF iteration is as follows.

\begin{enumerate}[label=\textbf{\arabic*.}]
  \item \textbf{Initialize.}  Draw $z_1,\ldots,z_N\sim\refd$ (fixed
        throughout) using \emph{sign-symmetric sampling}: draw $N/2$
        particles from $\refd|_{(0,L)}$ and mirror as
        $z_{N/2+j}=-z_j$.  This enforces $T_\theta(-z)=-T_\theta(z)$
        (bias gradients cancel by symmetry, keeping $g_\theta$ odd),
        so $\rho_\theta$ is symmetric about $x=0$.
        Initialize $\theta^0$ with small random weights (zero biases).

  \item \textbf{For $k=0,1,\ldots,K-1$:}
  \begin{enumerate}[label=(\alph*)]
    \item Integrate the augmented Neural ODE (Section~\ref{sec:alg:ode})
          for all $N$ particles with forward Euler ($N_{\rm ODE}$ substeps)
          to obtain
          $\{(x_i,\,\ell_i,\,J_i,\,\ell'_i)\}_{i=1}^N$.
    \item Evaluate $E(\theta^k) = \FQ+\FV+\FR$ by particle quadrature
          (Section~\ref{sec:alg:energy}) and compute its gradient
          $\nabla_\theta E$ by automatic differentiation through the
          forward ODE.
    \item Assemble $\Ghat(\theta^k)$ from the per-sample Jacobians
          $\partial_\theta T_\theta(z_i)$, cf.~\eqref{eq:G}.
    \item Solve $\Ghat(\theta^k)\,\xi = \nabla_\theta E(\theta^k)$
          by CG with backtracking (Section~\ref{sec:alg:solver}).
    \item Update $\theta^{k+1} = \theta^k
          - \alpha\,\min\!\bigl(1,\,C/\norm{\xi}\bigr)\,\xi$.
  \end{enumerate}

  \item \textbf{Reconstruct.}  Using the best parameters $\theta^*$
        (lowest energy seen), evaluate $T_{\theta^*}$ on a uniform grid,
        compute $\rhoT$, set $u=\sqrt{\rhoT}$, normalize in $L^2$,
        and export.
\end{enumerate}

\noindent
Since $E(\theta^k)$ is not guaranteed to decrease monotonically (the CG
solve is approximate and the step may overshoot), the best-parameter
tracking in step~3 is used for the final reconstruction.
\begin{remark}
    There are two alternatives in the algorithm. First, one may use MINRES instead of CG. It avoids using Tikhonov regularization. The stopping criteria in MINRES plays the regularization role. Second, there is no need to form $\hat{G}$ explicitly. The matrix vector multiplication can be directly implemented through $\partial_{\theta}T_{\theta}$.  
\end{remark}

\subsection{Product-map ansatz and dimension-specific choices}
\label{sec:alg:dims}

The algorithm above is dimension-independent.  For $d>1$ with a
separable potential, additive $V(\mathbf{x})=\sum_{k=1}^d V_k(x_k)$
or multiplicative $V(\mathbf{x})=\prod_{k=1}^d V_k(x_k)$,
we adopt a product-map ansatz:
\[
  T_\theta(\mathbf{z})
  = \bigl(T_{\theta_1}(z_1),\;\ldots,\;T_{\theta_d}(z_d)\bigr),
\]
where each $T_{\theta_k}$ is an independent copy of the 1D
boundary-preserving Neural ODE with its own parameters
$\theta_k\in\R^{M_1}$; the total parameter count is $M=d\cdot M_1$.
The metric tensor $\Ghat$ inherits a block-diagonal structure from the
product factorization, though our implementation assembles the full matrix.
In general, the product-map ansatz incurs a representation error since the
true ground-state density $\rho^*$ may not be a product.
We use this ansatz whenever the potential is separable, as in the tests of Sections~\ref{sec:results}--\ref{sec:3d}
and in~\ref{app:quadrature}, because it keeps the implementation simple:
each coordinate is an independent copy of the one-dimensional Neural ODE and
the density is the product $\rho_1(x_1)\cdots\rho_d(x_d)$, which is also
what makes it cheap to tabulate on a grid.  That product structure is
exactly its limitation.  When the potential couples the coordinates, the
ground state is not a product and the ansatz has an accuracy floor that no
parameter count can cross.  Sections~\ref{sec:nonsep}
and~\ref{sec:4d} measure that floor and replace the product map by a
full map $T_\theta\colon\Omega\to\Omega$ with coupled
components~\eqref{eq:fullmap}.

\paragraph{Reference densities}
The reference density $\refd$ must vanish at $\partial\Omega$
(Remark~\ref{rem:no-reg}), and the particle locations
$z_i\sim\refd$ should concentrate where the ground-state density is
large.  We use three choices, listed here for every test of
Section~\ref{sec:numerics}:

\begin{center}
\begin{tabular}{lll}
\toprule
Problem & Domain & Reference $\refd$ \\
\midrule
1D (\S\ref{sec:results}) & $(-1,1)$
   & Beta(2,2): $\tfrac{3}{4}(1-z^2)$ \\
2D (\S\ref{sec:2d_comparison}) & $(-16,16)^2$
   & Gaussian mixture per axis \\
3D (\S\ref{sec:3d}) & $(-8,8)^3$
   & Beta(5,5)$^3$: $C(1-z^2/L^2)^4$ \\
2D stirrer (\S\ref{sec:2d-stirrer})
   & $(-6,6)^2$
   & Beta(5,5) per axis \\
3D stirrer (\S\ref{sec:3d-stirrer})
   & $(-8,8)^3$
   & Beta(5,5) per axis \\
4D (\S\ref{sec:4d})
   & $(-6,6)^4$
   & Beta(5,5) per axis \\
\bottomrule
\end{tabular}
\end{center}

\noindent
For the 2D problem, $\refd_1(z) \propto (1-z^2/L^2)\sum_{k\in\mathcal{W}}
e^{-(z-k)^2/(2\sigma^2)}$ with
$\mathcal{W}=\{-12,-8,-4,0,4,8,12\}$, $\sigma=1.5$
(Section~\ref{sec:2d_comparison}).
The Beta(5,5) entries are the per-axis factor
of~\eqref{eq:beta55}, with $L$ the half-width of the domain and the
constant $C_1$ raised to the power $d$.

\subsection{Evaluating a trained map on a grid}
\label{sec:alg:eval}

{The output of PWGF is a map rather than an array, so using it as an initial
guess for a grid-based solver requires evaluating $\rho_\theta$ at prescribed
nodes.  For the product map of Section~\ref{sec:alg:dims} the density
factorizes,
$\rho_\theta(\mathbf{x})=\prod_{k=1}^{d}\refd_k(z_k)/T_{\theta_k}'(z_k)$ with
$z_k=T_{\theta_k}^{-1}(x_k)$, so it suffices to tabulate the $d$
one-dimensional profiles by integrating the map forward once; the cost of the
tabulation does not depend on the grid, and one table serves every grid.  A
map with coupled components (Section~\ref{sec:nonsep}) admits no such
factorization, and the density at a node is obtained by inverting the Neural
ODE: integrating the flow backwards from $\mathbf{x}$ returns
$\mathbf{z}=T_\theta^{-1}(\mathbf{x})$ and, accumulating
$\operatorname{tr}(\partial v/\partial w)$ along the same trajectory,
$\log|\det dT_\theta/d\mathbf{z}|$.  The density is then
$\rho_\theta(\mathbf{x})=\refd(\mathbf{z})e^{-\log|\det dT_\theta/d\mathbf{z}|}$.
Of the augmented state $(w,\ell,J,\ell')$ of
Section~\ref{sec:alg:ode}, only the flow position $w$ and the
log-determinant $\ell$ are carried here, $d+1$ scalars per node, since
evaluation needs the density but not the score.

The trained object is the discrete forward-Euler map, so the inversion must
use the implicit step $w_j = w_{j+1}-\Delta s\,v(w_j)$, solved by fixed-point
iteration; an explicit backward step inverts a different map and misplaces the
density by several percent.
}

\section{Numerical Results}
\label{sec:numerics}

In 2D, 3D, and 4D experiments below, we first run PWGF to obtain an approximate
ground state, and then use it as a warm start for the $H^1$ Sobolev
gradient flow~\cite{zhang2024sobolev}. Thus we refer to the PWGF output as the
\emph{PWGF warm start} throughout.
The three tests of
Sections~\ref{sec:results}--\ref{sec:3d} use the exact solution in 1D and
finite difference reference solutions in 2D and 3D, and both the PWGF
training and the $H^1$ flow run on the CPU of a MacBook Pro with an Apple M1
processor; their hyperparameters are collected in
Table~\ref{tab:hyperparams}.  The PWGF training in these
three tests runs in single precision; \ref{app:quadrature} examines
double-precision and Gauss--Legendre quadrature variants.  All three use the same learning rate $\alpha=0.005$; the 1D and
2D experiments use clip $C=10$ while the 3D experiment uses $C=50$.
The 2D problem additionally uses 4-fold sign-symmetric sampling
(Section~\ref{sec:2d_comparison}) to enforce the symmetry of the potential and
keep network biases at zero throughout optimization.

The larger tests of Sections~\ref{sec:nonsep}
and~\ref{sec:4d}, and the $999^3$ comparison of
Section~\ref{sec:3d-warmstart}, instead run on one NVIDIA GH200 GPU, with a
wider network and a larger training budget stated there.
The reference solutions they are measured against use $Q^{20}$
spectral elements in three dimensions, $Q^{10}$ in four, and an $800^2$ finite
difference grid in two.

\begin{table}[ht]
\centering
\caption{Hyperparameters for the three experiments of
Sections~\ref{sec:results}--\ref{sec:3d}.  $M_1=H^2+4H+1$
(parameter count for a $1{\to}H{\to}H{\to}1$ network);
$M=d\cdot M_1$ for $d$-dimensional product maps.}\label{tab:hyperparams}
\begin{tabular}{lccc}
\toprule
Parameter & 1D & 2D & 3D \\
\midrule
$N$ (particles)       & 3000  & 3000  & 6000  \\
$K$ (steps)           & 400   & 400   & 400   \\
$\alpha$ (step size)  & 0.005 & 0.005 & 0.005 \\
$H$ (hidden width)    & 10    & 10    & 10    \\
$M$ (parameters)      & 141   & 282   & 423   \\
$N_{\rm ODE}$ (Euler steps) & 10    & 10    & 10    \\
$n_{\rm CG}$ (CG iterations)  & 100   & 200   & 300   \\
$E_{\rm tol}$ (backtracking) & 0.05  & 0.05  & 5.0   \\
Clip $C$              & 10    & 10    & 50    \\
\bottomrule
\end{tabular}
\end{table}

\subsection{1D Test}
\label{sec:results}

We solve the stationary GPE on $\Omega=(-1,1)$ with Dirichlet boundary
conditions:
\begin{equation}
  -u'' + V(x)\,u + \beta|u|^2 u = \lambda u,
  \quad x\in(-1,1),\quad u(\pm 1)=0,\quad \norm{u}_{L^2}=1.
  \label{eq:gpe_eigen}
\end{equation}
The potential and interaction coefficient are
\[
  V(x) = \beta\cos^2\!\Bigl(\frac{\pi(x+1)}{2}\Bigr),\qquad \beta=10.
\]
This problem has the exact ground-state solution
\begin{equation}
  u^*(x) = \sin\!\Bigl(\frac{\pi(x+1)}{2}\Bigr),\qquad
  \lambda^* = \frac{\pi^2}{4} + \beta,\qquad
  E^* = \frac{\lambda^*}{2} - \frac{3\beta}{16}.
  \label{eq:exact}
\end{equation}
Numerically: $\lambda^* \approx 12.4674$, $E^* \approx 4.3587$.

\begin{figure}[!ht]
  \centering
  \includegraphics[width=\textwidth]{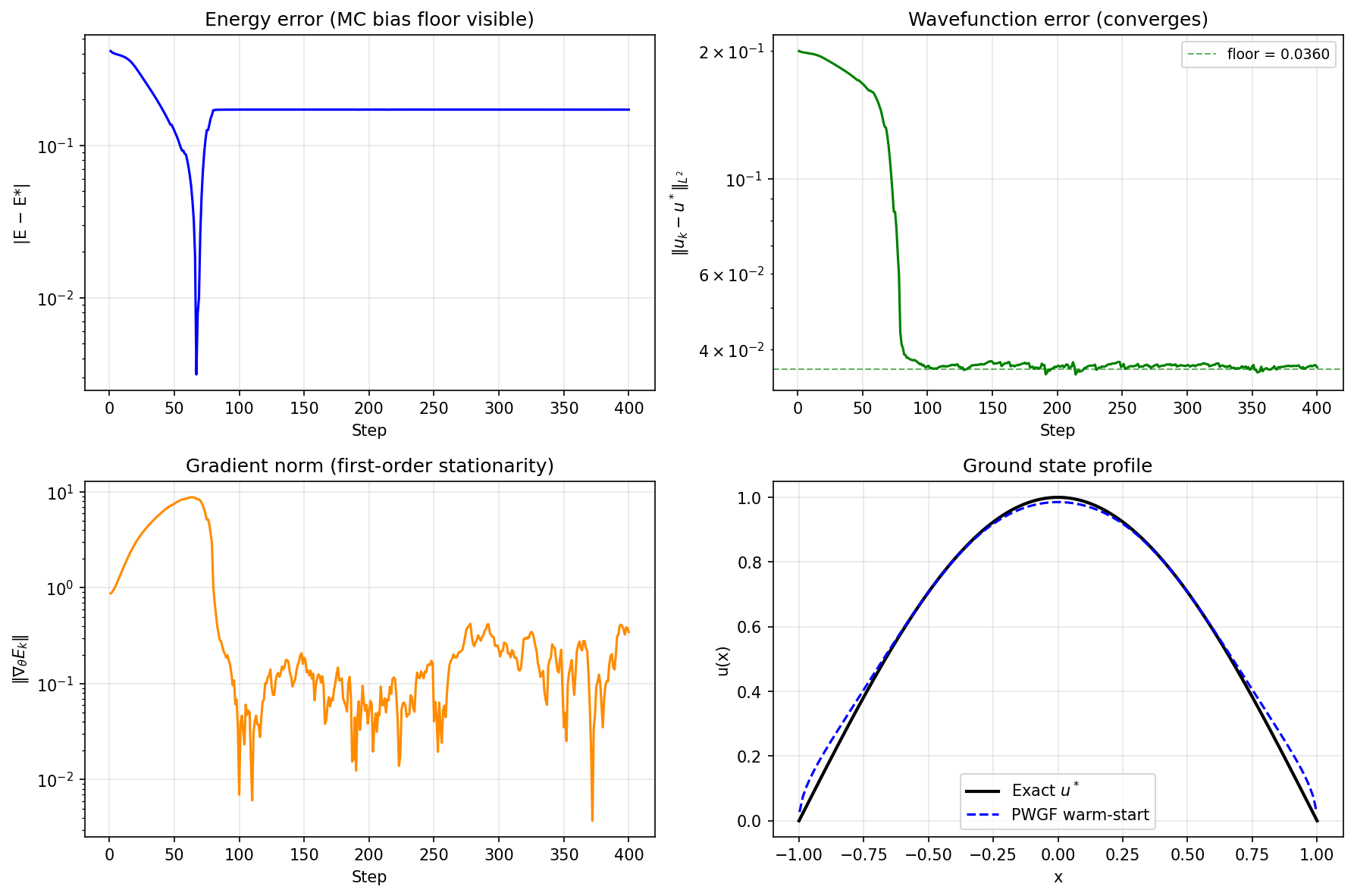}
  \caption{1D PWGF ground-state computation.
    \emph{Top-left:} energy error $|E^k-E^*|$ (log scale).
    \emph{Top-right:} $L^2$ error $\|u^k-u^*\|$; decreases to the
    floor $\approx 0.036$ (Table~\ref{tab:errorsummary}).
    \emph{Bottom-left:} gradient norm $\|\nabla_\theta E^k\|$.
    \emph{Bottom-right:} reconstructed $u_\theta$  vs.\
    exact $u^*$.}
  \label{fig:results}
\end{figure}

Figure~\ref{fig:results} shows the PWGF iterates with the settings in
Table~\ref{tab:hyperparams} (runtime $\approx 26$\,s on MacBook Pro with M1 CPU).
The energy estimate plateaus at $E\approx 4.187$, below $E^*\approx 4.359$,
because the particle-based Fisher estimator systematically underestimates
$F_Q$ due to tail under-sampling (see Section~\ref{sec:ablation}
and \ref{sec:fisher1d}).
The wavefunction error $\|u^k-u^*\|_{L^2}$ decreases from $0.20$ at step~1
to $0.036$ by step~100 and remains stable thereafter. 

\begin{remark}[Energy below $E^*$]
The particle estimate $E_{\rm PWGF}\approx 4.19$ falls below the exact
ground-state energy $E^*\approx 4.36$.  This does not violate the variational
principle: the true energy $E(\rho_\theta)$ is above $E^*$, but its Monte
Carlo estimate is biased downward.  All three energy components $\FQ$, $\FV$,
$\FR$ are estimated by particle quadrature and are subject to bias; the
ablation study in Section~\ref{sec:ablation} identifies $\FQ$ (Fisher
information) as the dominant source.
\ref{sec:fisher1d} sharpens this: with the Beta(2,2)
reference the true $E(\rho_\theta)$ is not merely above $E^*$ but
\emph{infinite}, since the Fisher integral of every representable density
diverges logarithmically at the boundary, so any finite particle estimate
necessarily sits below it.
\end{remark}

\subsubsection*{Sensitivity to hyperparameters}
\label{sec:ablation}

\begin{figure}[!ht]
  \centering
  \includegraphics[width=\textwidth]{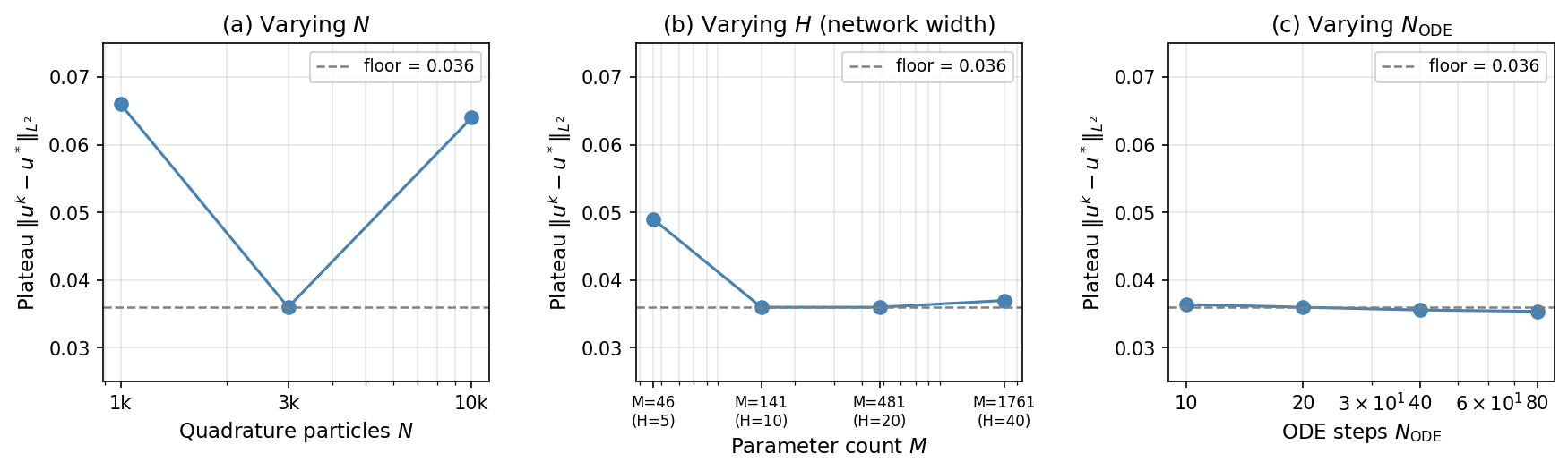}
  \caption{Plateau $\|u^k-u^*\|_{L^2}$ vs.\ each hyperparameter
    (baseline: $N=3000$, $H=10$, $N_{\rm ODE}=10$; dashed line at $0.036$).
    \emph{(a)}~$N$: non-monotone; $N=3000$ is optimal.
    \emph{(b)}~$H$: saturates at $H=10$.
    \emph{(c)}~$N_{\rm ODE}$: negligible improvement
    ($0.0364\to0.0354$ for $8\times$ more compute).}
  \label{fig:ablation}
\end{figure}

To identify the dominant source of the $\approx 3.6\%$ error floor,
we vary each hyperparameter independently while holding the others fixed
(Figure~\ref{fig:ablation} and Table~\ref{tab:errorsummary}).

\begin{table}[!ht]
\centering
\begin{tabular}{llll}
\toprule
Axis varied & Range & Plateau range & Observations \\
\midrule
$N$ (particles)     & $1000\to10000$ & $0.036$--$0.066$ & non-monotone  \\
$H$ (network width) & $5\to40$       & $0.036$--$0.049$ & improves $H{=}5\to10$ \\
$N_{\rm ODE}$ (ODE steps) & $10\to80$ & $0.0354$--$0.0364$ & little improvement \\
\bottomrule
\end{tabular}
\caption{Error-floor ablations (Figure~\ref{fig:ablation}).
  The bottleneck is MC bias in $\nabla_\theta F_Q$ from tail
  under-sampling.}
\label{tab:errorsummary}
\end{table}

All three ablations identify the same bottleneck: systematic MC bias in
the Fisher estimator $F_Q$ from tail under-sampling.  Increasing $N$
beyond $3000$ does not help because the optimizer converges more tightly
to the biased MC minimum; network capacity saturates at $H=10$; and
$8\times$ more ODE steps yield only $2.7\%$ improvement.
See~\ref{sec:fisher1d} for how Gauss--Legendre quadrature
and a better reference density improve the performance.

\subsection{2D Test}
\label{sec:2d_comparison}

The second test problem is the GPE on $(-16,16)^2$ with $\beta=10$:
\begin{equation}
  -\Delta u + V(x_1,x_2)\,u + \beta|u|^2 u = \lambda u,
  \quad (x_1,x_2)\in(-16,16)^2,\quad \beta=10,
  \label{eq:gpe2d_comp}
\end{equation}
with potential
$V(x_1,x_2) = 2\sin^2(\pi x_1/4)\sin^2(\pi x_2/4)$.
No analytical ground state is known; the reference value
$E^*\approx0.2171$ is computed by the Sobolev gradient flow method
of~\cite{zhang2024sobolev}.

\subsubsection{Reference density}

The potential vanishes on the lines
$\{x_1\in4\mathbb{Z}\}\cup\{x_2\in4\mathbb{Z}\}$; the ground-state density
concentrates at their $7\times7$ interior intersections
$(x_1,x_2)=(4j,4k)$, $j,k\in\{-3,\ldots,3\}$, forming 49 bumps.
We use a \emph{Gaussian-mixture reference density} that already places mass at
all seven potential wells per coordinate:
\begin{equation}
  \refd_1(z)
  = C\cdot\Bigl(1-\frac{z^2}{L^2}\Bigr)
    \sum_{k\in\mathcal{W}}\exp\!\Bigl(-\frac{(z-k)^2}{2\sigma^2}\Bigr),
  \quad
  \mathcal{W}=\{-12,-8,-4,0,4,8,12\}, 
  \label{eq:gmix_ref}
\end{equation}
where $\sigma=1.5$, the boundary factor enforces $\refd_1(\pm L)=0$ and $C$ is a
normalization constant.  The 2D reference is
$\refd(\mathbf{z})=\refd_1(z_1)\,\refd_1(z_2)$, and its score
\begin{equation}
  \partial_z\log\refd_1(z)
  = \sum_{k\in\mathcal{W}}w_k(z)\cdot\frac{-(z-k)}{\sigma^2}
    - \frac{2z}{L^2-z^2},
  \label{eq:gmix_score}
\end{equation}
where $w_k(z)$ are softmax-normalized Gaussian weights, enters the augmented
ODE for the Fisher estimator.

The symmetry $V(-x_1,x_2)=V(x_1,-x_2)=V(x_1,x_2)$ implies the optimal maps
are odd: $T_k(-z)=-T_k(z)$.  We enforce this analytically via
\emph{4-fold sign-symmetric sampling}: draw $N/4$ base particles
$(z_1^+,z_2^+)$ from the positive quadrant and include all four sign mirrors
$(\pm z_1^+,\pm z_2^+)$.  Gradient contributions of each bias parameter cancel
identically across opposite pairs, keeping the network biases at zero
throughout the iteration.

\subsubsection{Results}

With the settings in Table~\ref{tab:hyperparams}, the energy decreases from
$\approx0.225$ at step~20 to $E_{\rm best}=0.2162$ by step~400, with the
CG step accepted at every iteration.  The $7\times7$ multi-bump structure
is reproduced (Figure~\ref{fig:warmstart_combined}). CPU runtime is $\,78$\,s on MacBook Pro with M1 CPU.

\subsubsection{Warm-start comparison with the H\texorpdfstring{$^1$}{1} gradient flow}
\label{sec:h1flow}

The PWGF output is bilinearly interpolated onto the $200^2$ FD grid and
used to initialize the $H^1$ gradient flow:

\begin{center}
\begin{tabular}{lll}
\toprule
Initialization & Initial FD energy $E^{(0)}$ & $|E^1 - E^*_h|$ \\
\midrule
Constant one   & $0.6495$  & $8\times10^{-3}$ \\
Random ($|\mathcal{N}(0,1)|$, seed 42) & $29.49$ & $1.5\times10^{-2}$ \\
PWGF warm start & $\mathbf{0.2327}$ & $\mathbf{1.3\times10^{-3}}$ \\
\bottomrule
\end{tabular}
\end{center}

\begin{figure}[!ht]
  \centering
  \includegraphics[width=\textwidth]{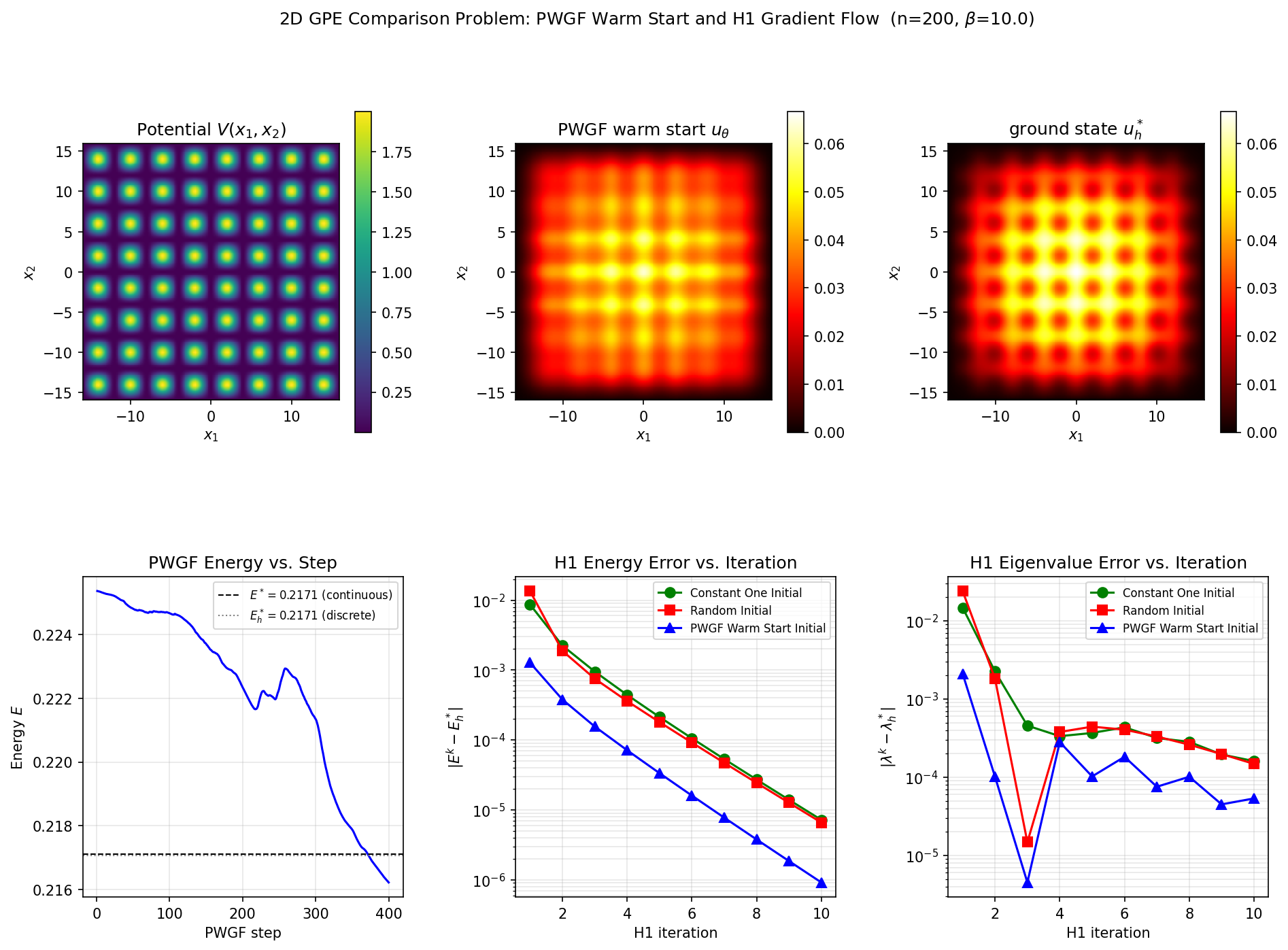}
  \caption{2D problem ($\beta=10$, $200^2$).
    \emph{Top:} potential $V$, PWGF warm start $u_\theta$
    and the ground state $u^*_h$.
    \emph{Bottom left:} PWGF energy history.
    \emph{Bottom center/right:} energy and eigenvalue error vs.\  $H^1$
    iteration.}
  \label{fig:warmstart_combined}
\end{figure}

At H1 step~1, the PWGF warm start is ${\sim}7\times$ closer to $E^*_h$
than constant-one and ${\sim}11\times$ closer than random
(Figure~\ref{fig:warmstart_combined}).

\ref{sec:2d-quadrature} repeats this computation
with the Monte Carlo estimates replaced by per-axis Gauss--Legendre
quadrature.  The estimator bias disappears, but the accuracy is limited by
the product-map representation (the rank-one floor of $u^*_h$ is $0.048$)
and the $H^1$ warm-start value is unchanged, which supports retaining the
setup of Table~\ref{tab:hyperparams}.

\subsection{3D Test}
\label{sec:3d}

For the third test we take the GPE on $(-8,8)^3$ with $\beta=1600$:
\begin{equation}
  -\Delta u + V(x_1,x_2,x_3)\,u + \beta|u|^2 u = \lambda u,
  \quad \mathbf{x}\in(-8,8)^3,\quad \beta=1600,
  \label{eq:gpe3d}
\end{equation}
with harmonic trap plus optical lattice
\begin{equation}
  V(x_1,x_2,x_3)
  = x_1^2+x_2^2+x_3^2
  + 100\!\Bigl(\sin^2\tfrac{\pi x_1}{4}
               +\sin^2\tfrac{\pi x_2}{4}
               +\sin^2\tfrac{\pi x_3}{4}\Bigr).
  \label{eq:V3d}
\end{equation}
Reference values from \cite{zhang2024sobolev} are $E^*=33.80228$ and $\lambda^*=80.895$.

\subsubsection{Reference density}

We use a scaled Beta(5,5) reference density on each coordinate interval $(-L,L)$,
concentrated within $|z_k|<4$ where the potential is low:
\begin{equation}
  \refd(\mathbf{z})
  = C_1^3\prod_{k=1}^{3}\Bigl(1-\frac{z_k^2}{L^2}\Bigr)^{\!4},
  \quad C_1 = \frac{315}{256L},
  \label{eq:beta55}
\end{equation}
with score $\partial_{z_k}\log\refd = -8z_k/(L^2-z_k^2)$.
The normalization $C_1=315/(256L)$ follows from $\int_{-1}^1(1-t^2)^4\,dt=256/315$,
obtained by expanding $(1-t^2)^4$ and integrating term by term.
Quadrature uses $N=6000$ particles drawn with 8-fold sign symmetry
($750$ base triplets $\times$ all $(\pm,\pm,\pm)$ sign patterns), preserving
the sign-symmetry of $\nabla_\theta E$.
All three coordinate networks are initialized identically with small random
weights (scale $0.01$) to avoid the saddle point at zero initialization.

\subsubsection{Results}

\begin{figure}[!ht]
  \centering
  \includegraphics[width=\textwidth]{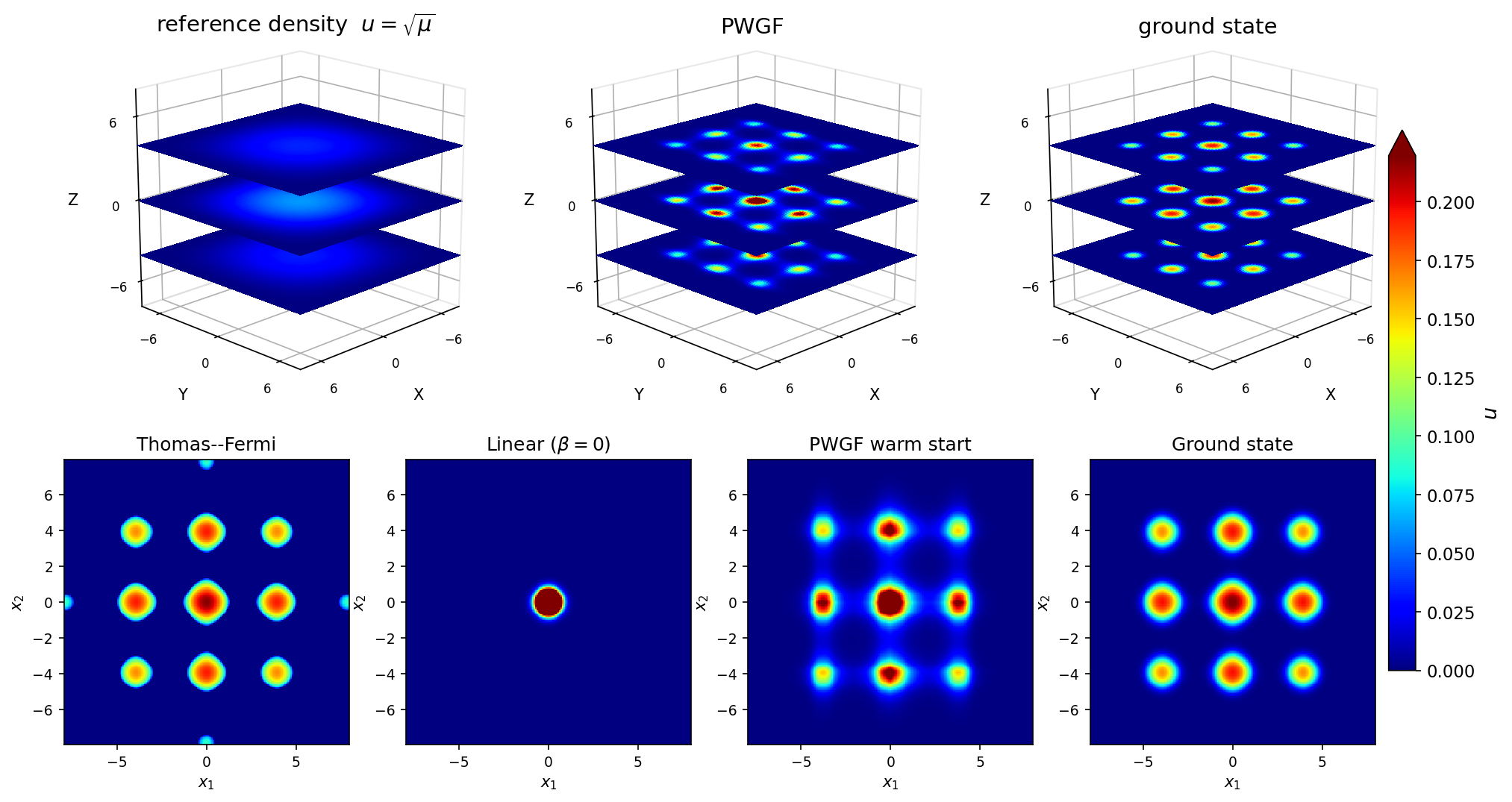}
  \caption{{Both rows share one colour scale, capped at the peak
    of the ground state; the linear state, whose peak is $1.97$, and the PWGF
    output run past the cap.
    \emph{Top row:} slices at $z\approx-4,0,+4$.
    \textbf{Left:} the Beta(5,5)$^3$ reference density.
    \textbf{Middle:} the PWGF output.  \textbf{Right:} the
    ground state, $3\times3\times3$ separated blobs.
    \emph{Bottom row:} the $z=0$ midplane of each initialization and of the
    ground state.  The Thomas--Fermi and linear
    initializations in the bottom row are defined in
    Section~\ref{sec:3d-warmstart}.}}
  \label{fig:3d-slices}
  \label{fig:3d-visualization-n199}
  \label{fig:3d-reference-density}
  \end{figure}

With the settings in Table~\ref{tab:hyperparams}, the energy decreases from
$E\approx84$ at step~20 to $E_{\rm best}\approx38.70$
(within $15\%$ of $E^*$)
by step~200.
The ground state has a $3\times3\times3$ blob structure dictated by the
optical lattice (Figure~\ref{fig:3d-visualization-n199}).  The PWGF solution after $K=400$ steps shown in Figure~\ref{fig:3d-slices} is certainly still far from the true ground state, but it is much better than 
unimodal Beta(5,5) shown in Figure~
\ref{fig:3d-reference-density}. In other words, this 3D test does show that PWGF can recover some interesting structure about the true minimizer of GP energy, though it is quite difficult to use PWGF to obtain a highly accurate solution.

Despite this qualitative mismatch, the PWGF warm start still provides a
useful initialization for the H1 flow because the total mass is placed in
the correct spatial region ($|x_k|<4$), even though its fine structure is
wrong.

\subsubsection{Warm-start comparison with the H\texorpdfstring{$^1$}{1} gradient flow}
\label{sec:3d-warmstart}

The PWGF output (on a $40^3$ grid) is trilinearly interpolated onto the
$99^3$ FD grid and used to initialize the H1 gradient flow.

\begin{figure}[!ht]
  \centering
  \includegraphics[width=\textwidth]{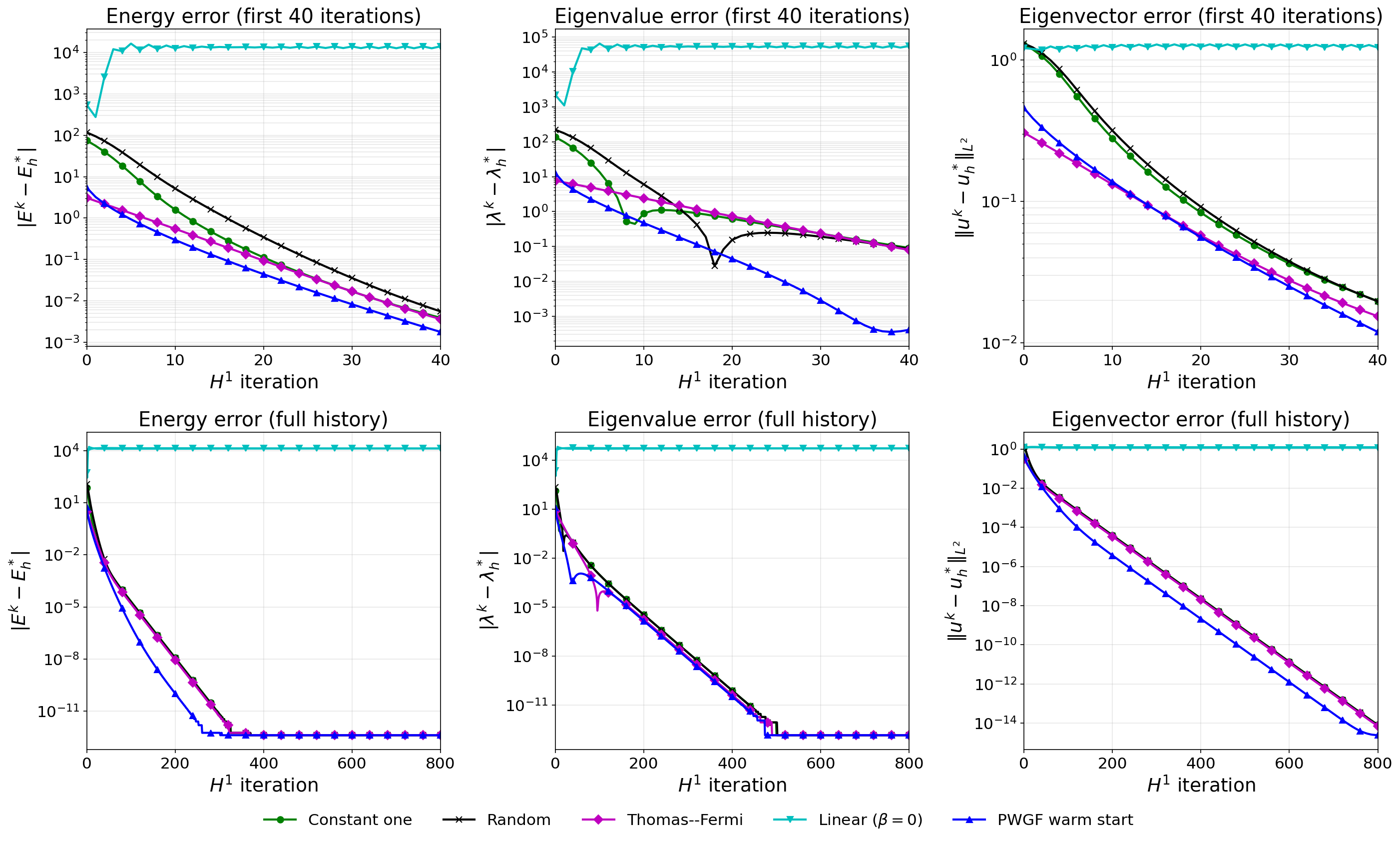}
  \caption{{3D GPE ($\beta=1600$) on the $99^3$ finite-difference
    grid.  \emph{Top row:} the first $40$ iterations.
    \emph{Bottom row:} the same runs over the full $800$ iterations, where they
    have merged onto a common floor.  Errors are measured
    against the minimizer of the \emph{discrete} energy on this grid,
    $E^*_h=33.769419797696$ and $\lambda^*_h=80.855026128719$, obtained by
    running the same flow from the Thomas--Fermi start for $1600$ iterations to
    a residual of $5.6\times10^{-14}$.}}
  \label{fig:3d-warmstart-combined}
\end{figure}

The PWGF warm start enters the H1 flow ${\sim}17\times$
closer to $E^*$
than the constant-one cold start and ${\sim}29\times$ closer
than the random
cold start at step~1.  After 10 H1 steps, the warm start reaches
$E=34.07$
(within $1\%$ of $E^*$), while constant-one reaches
$35.36$ and random
reaches $38.93$ (Figure~\ref{fig:3d-warmstart-combined}).

The PWGF phase runs in ${\approx}\,240$\,s (MacBook Pro with M1 CPU, $K{=}400$ steps,
$M{=}423$ parameters) for the 3D problem. 

\begin{table}[!ht]
\centering
\caption{Modified $H^1$ iterations required to reach a given
energy error on the $99^3$ grid, measured against the
minimizer $E^*_h$ of the discrete energy quoted in
Figure~\ref{fig:3d-warmstart-combined}, i.e.\ the horizontal offset of
Figure~\ref{fig:3d-warmstart-combined} made quantitative.
The modified $H^1$ flow carries the metric
$(\nabla w,\nabla z)+a(w,z)$ with $a=20$ and $dt=0.1$.}
\label{tab:h1-iters-n99}
{\small
\begin{tabular}{lrrrrr}
\toprule
Initialization & $E^{(0)}$ & $10^{-2}$ & $10^{-4}$ & $10^{-6}$ & $10^{-8}$\\
\midrule
Constant one       & $108.40$ & $34$ & $80$ & $141$ & $203$\\
Random             & $150.92$ & $37$ & $79$ & $140$ & $202$\\
Thomas--Fermi      & $36.84$  & $34$ & $76$ & $137$ & $199$\\
Linear ($\beta=0$) & $572.71$ & --- & --- & --- & ---\\
PWGF warm start    & $\mathbf{39.19}$ & $\mathbf{29}$ & $\mathbf{61}$ & $\mathbf{99}$ & $\mathbf{144}$\\
\bottomrule
\end{tabular}
}
\end{table}

\ref{sec:3d-quadrature} examines how far the
product ansatz itself can be pushed on this problem.  

\begin{figure}[!ht]
  \centering
  \includegraphics[width=\textwidth]{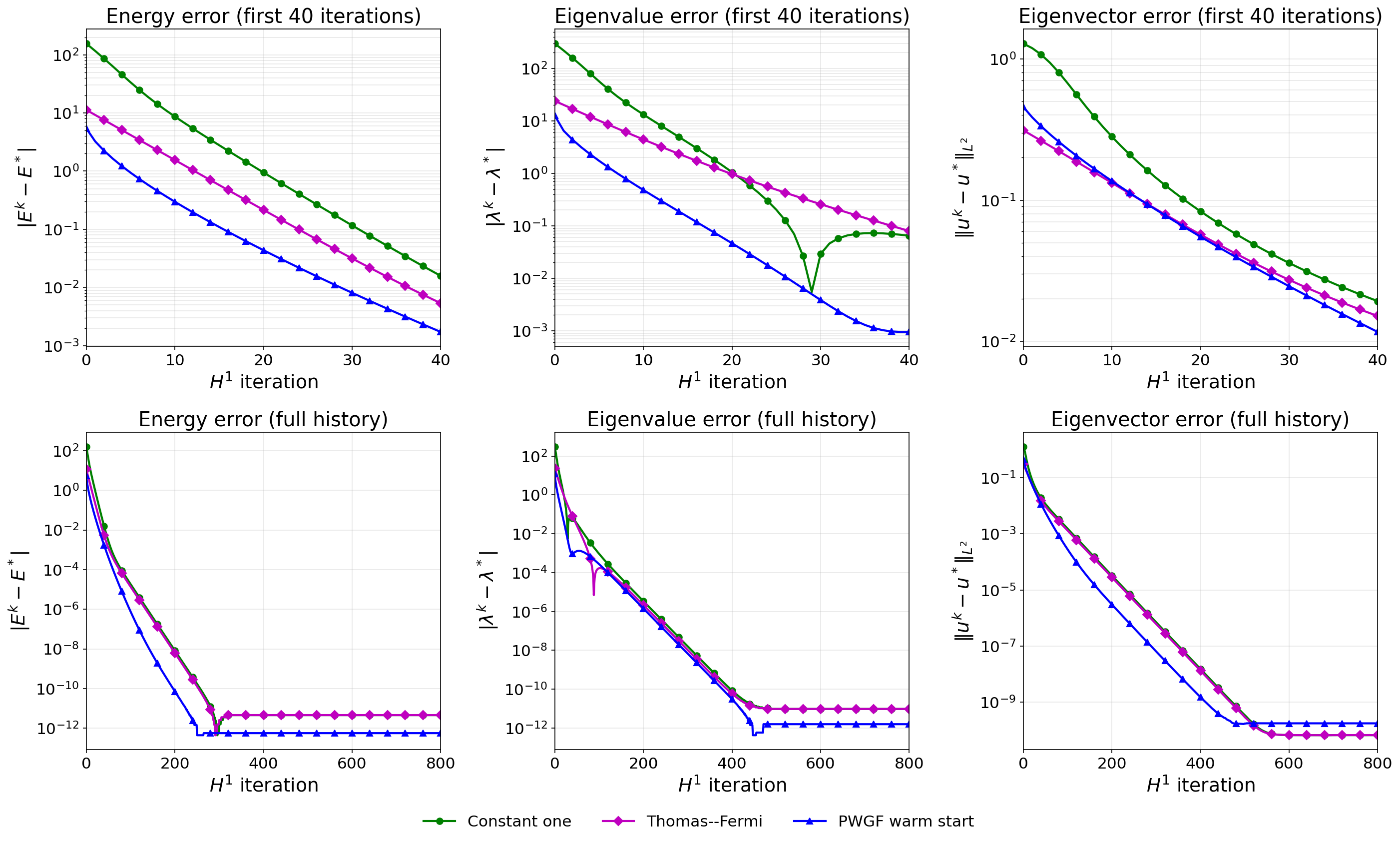}
  \caption{Complete modified $H^1$ history on the $999^3$
    grid ($a=20$, $dt=0.1$), double precision,
    against~\eqref{eq:ongrid-ref}.  \emph{Top row:} the first $40$
    iterations.  \emph{Bottom row:} $800$ iterations, and the horizontal offset is the saving
    as quantified in Table~\ref{tab:h1-iters-n999}.  The linear start does not converge and
    is omitted.  }
  \label{fig:3d-cheap-inits}
\end{figure}

{\medskip\noindent\textbf{Cheap standard initializations.}
Constant-one and random data are the weakest baselines, so we now compare
four initializations on a $999^3$ grid with about $10^9$ degrees of
freedom, fine enough that the cost of one mesh-based iteration is
representative of large-scale practice.  We use the $Q^{20}$ spectral element
discretization of~\cite{li2022sem,zhang2024sobolev} and the GPU tensor-product
solver of~\cite{liu2024gpucicp,liu2025gpu}, and run the modified $H^1$ flow (metric
$(\nabla w,\nabla z)+a(w,z)$, $dt=0.1$, $a=20$).  The separable solves are
one-dimensional eigendecompositions along each axis, so no $d$-dimensional
operator is assembled, and the linear ($\beta=0$) initialization is computed
by the preconditioned inverse iteration of~\cite{liu2025gpu}.  All runs are in
double precision so that the convergence plateau is not masked by
arithmetic.

\medskip\noindent\emph{Reference values.}  Errors are measured against the minimizer
of the discrete energy on this grid, obtained by running the flow in double
precision to residual stall:
\begin{equation}
  E^* = 33.802279005474, \qquad \lambda^* = 80.895114406041 .
  \label{eq:ongrid-ref}
\end{equation}
These agree with the converged values of~\cite{zhang2024sobolev} to
$1.4\times10^{-12}$ and $1.1\times10^{-11}$.

\medskip\noindent\emph{Thomas--Fermi (TF).}  In the strongly interacting
regime the kinetic term is negligible outside a thin boundary layer.
Dropping $-\Delta u$ from~\eqref{eq:gpe3d} leaves the algebraic equation
$(V+\beta u^2)u=\lambda u$, whose nontrivial branch is
\begin{equation}
  u_{\rm TF}(\mathbf{x})
  = \sqrt{\frac{\bigl(\lambda_{\rm TF}-V(\mathbf{x})\bigr)_+}{\beta}},
  \qquad (s)_+=\max(s,0),
  \label{eq:tf}
\end{equation}
with $\lambda_{\rm TF}$ the unique root of the increasing function
$\Phi(\lambda)=\beta^{-1}\!\int_\Omega(\lambda-V)_+\,d\mathbf{x}-1$.  No PDE
is solved.  The construction costs $1.6$\,s on $10^9$ points and gives
$\lambda_{\rm TF}=75.97$.  This is the standard large-$\beta$ initialization
for BEC computations~\cite{bao2004,baocai2013}.
Its energy depends on the grid.  The profile has a square-root
free boundary where $V=\lambda_{\rm TF}$, at which $|\nabla u_{\rm TF}|^2$ is
not integrable, so the discrete Dirichlet energy grows logarithmically as the
mesh resolves that edge: on the finite difference grids it rises from $36.84$
at $n=99$, the value in Table~\ref{tab:h1-iters-n99}, to $39.81$ at $n=399$,
and the $Q^{20}$ elements used here resolve the edge more sharply still, which
is why $E^{(0)}$ for TF is $45.08$ on this grid.  The remaining
initializations are smooth and show no such drift.

\medskip\noindent\emph{Linear ($\beta=0$) ground state.}  The lowest eigenpair
of $-\Delta+V$ by preconditioned inverse iteration.  Since $V$ in~\eqref{eq:V3d} is
additively separable, each solve is a tensor-product operation.  The
construction costs $54$\,s, or $39$\,s at the looser PCG tolerance that
suffices for an initial guess; it is the only initialization here that solves
a PDE to build itself, and the only setup that costs more than a couple of
seconds.

\medskip\noindent\emph{PWGF warm start.}  The product-map ansatz factorizes
the PWGF density exactly,
\begin{equation}
  \rho_\theta(x_1,x_2,x_3)=\prod_{k=1}^{3}\frac{\refd_1(z_k)}{T_k'(z_k)},
  \qquad z_k=T_k^{-1}(x_k),
  \label{eq:pwgf-rank1}
\end{equation}
so $u_\theta$ is a rank-one product of three 1D profiles, which we tabulate
at $4001$ points per axis and evaluate on whatever grid is required.
Evaluating the map gives initial energies of $39.19$, $39.21$ and $39.23$
on the $99^3$ and $199^3$ finite difference grids
and on this grid, a spread under $0.1\%$, as a mesh-free function must give.
The evaluation on $10^9$ points costs $16$\,ms, two orders of magnitude below
the cheapest competing construction.

\medskip
Table~\ref{tab:3d-cheap-inits} reports the cost and initial quality of each
start and Figure~\ref{fig:3d-cheap-inits} the errors along the flow.  The
linear start is unusable at $\beta=1600$, its energy growing to
$5\times10^{3}$, and is omitted from the figure.  Of the rest PWGF begins
closest, $E^{(0)}=39.23$ against $45.08$ for TF, and keeps the smallest energy
error at every one of the $40$ iterations, $1.7\times10^{-3}$ at $k=40$
against $5.4\times10^{-3}$ for TF and $1.6\times10^{-2}$ for constant-one.
Taking energy, eigenvalue and eigenvector together, PWGF converges faster
than TF: it overtakes TF in $L^2$ at iteration $14$ and stays
ahead for the rest of the first $40$ iterations (Figure~\ref{fig:3d-cheap-inits}) and
reaches $10^{-8}$ in $142$ iterations against $194$
(Table~\ref{tab:h1-iters-n999}).  On cost, one $H^1$ iteration costs
$2.50$\,s in double precision on one NVIDIA GH200 GPU with our Python
JAX implementation, so the $52$ iterations saved against TF are $130$\,s of GPU
time.  They are bought
with the $240$\,s of offline training, which runs on the CPU of a MacBook Pro
with an Apple M1 processor and never touches the grid, its cost being set by
$N$ and $M$ instead.  Counted as bare seconds the warm start returns a little
over half its own cost at $999^3$, and that accounting is the conservative
one, since it charges laptop CPU seconds against GH200 GPU seconds and
prices a map
that is trained once against a grid on which it is evaluated in $16$\,ms.  The rank-one structure that makes it cheap
to evaluate is also what bounds it, and removing the axis-aligned ridges
requires the full map of Section~\ref{sec:nonsep}.
}

\begin{table}[!ht]
\centering
\caption{3D GPE ($\beta=1600$) on the $999^3$ $Q^{20}$ SEM
grid ($\approx10^9$ DoFs, double precision, one NVIDIA GH200 GPU): cost and
initial
quality of each initialization against the ground state reference
solution~\eqref{eq:ongrid-ref}.
Setup time excludes the $240$\,s
offline PWGF optimization.  ``Grid-dependent'' records
whether the construction must be redone when the grid changes; the PWGF map
is trained once, independently of any grid, and only its evaluation
($16$\,ms) touches the grid.}
\label{tab:3d-cheap-inits}
{\begin{tabular}{lrrrr}
\toprule
Initialization & Setup (s) & Grid-dependent? & $E^{(0)}$ & Gap above $E^*$\\
\midrule
Constant one         & $0.01$ & yes & $189.33$ & $460\%$\\
Thomas--Fermi        & $1.6$  & yes & $45.08$  & $33.4\%$\\
Linear ($\beta=0$)   & $54$   & yes & $548.33$ & $1522\%$\\
PWGF warm start      & $\mathbf{0.02}$ & no  & $\mathbf{39.23}$ & $\mathbf{16.0\%}$\\
\bottomrule
\end{tabular}
}
\end{table}

\begin{table}[!ht]
\centering
\caption{Modified $H^1$ iterations to reach a given energy
error on the $999^3$ grid.  The pattern matches the $99^3$ grid
(Table~\ref{tab:h1-iters-n99}), so the saving is essentially
grid-independent: PWGF reaches $10^{-8}$ in $142$ iterations against $194$
and $198$, a saving of $52$--$56$ iterations, about $27\%$. }
\label{tab:h1-iters-n999}
{\small
\begin{tabular}{lrrrrr}
\toprule
Initialization & $E^{(0)}$ & $10^{-2}$ & $10^{-4}$ & $10^{-6}$ & $10^{-8}$\\
\midrule
Constant one    & $189.33$ & $43$ & $79$ & $138$ & $198$\\
Thomas--Fermi   & $45.08$  & $37$ & $75$ & $134$ & $194$\\
PWGF warm start & $\mathbf{39.23}$ & $\mathbf{29}$ & $\mathbf{61}$ & $\mathbf{98}$ & $\mathbf{142}$\\
\bottomrule
\end{tabular}
}
\end{table}

\subsection{Non-separable potentials: product map versus full map}
\label{sec:nonsep}

{Every test so far has used the product-map ansatz of
Section~\ref{sec:alg:dims}, for which $\hatG$ is block diagonal and the
representable densities are the rank-one products $\rho_1(x_1)\rho_2(x_2)$.
This section considers test cases in which neither holds: the potentials are
not separable, and we compare the product map against a full map
$T_\theta\colon\Omega\to\Omega$ with coupled components and dense $\hatG$,
in two dimensions and then at $999^3$ in three.
}

\subsubsection{2D stirrer potential}
\label{sec:2d-stirrer}

\begin{figure}[!ht]
  \centering
  \includegraphics[width=\textwidth]{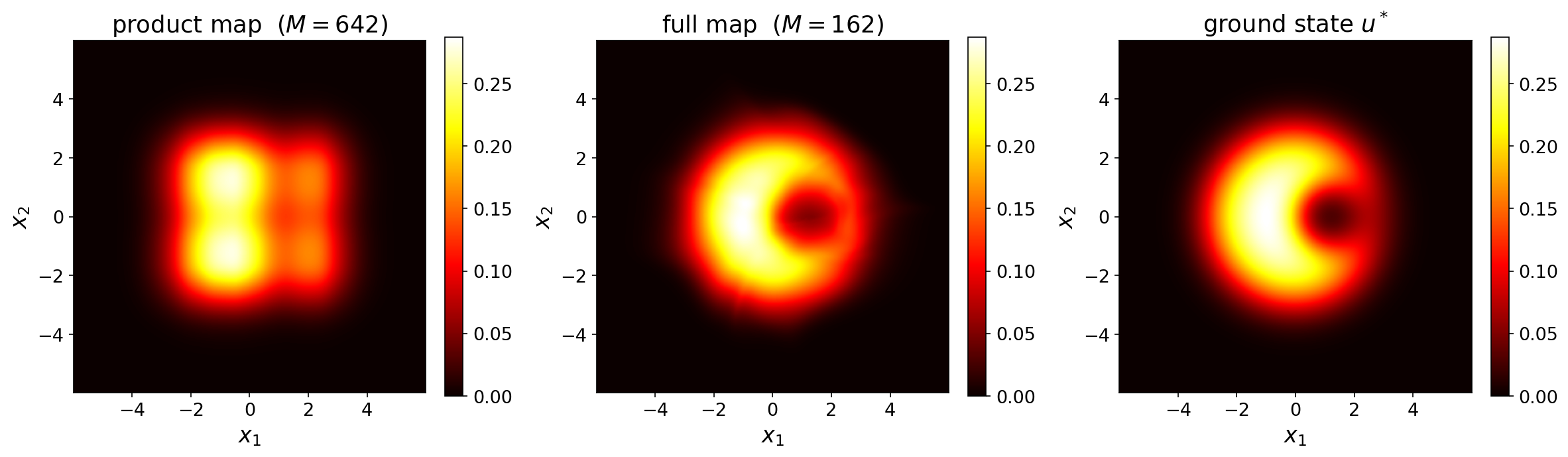}
  \caption{2D stirrer potential: the product-map density at $M=642$,
    the full-map density at $M=162$ and the ground state $u^*$.  }
  \label{fig:stirrer2d}
\end{figure}

{\noindent\emph{Problem.}  We use the harmonic-plus-stirrer potential
of~\cite{AntoineDuboscq2015,bao2004}, a harmonic trap perturbed by a localized
Gaussian obstacle:
\begin{equation}
  V(x_1,x_2) = x_1^2 + x_2^2
  + 2w_0\,e^{-\delta\left((x_1-r_0)^2 + x_2^2\right)},
  \label{eq:stirrer2d}
\end{equation}
with $w_0=10$, $\delta=1$, $r_0=1$, $\beta=100$ on $\Omega=(-6,6)^2$.  The
obstacle couples the two coordinates, so $V$ cannot be written as
$V_1(x_1)+V_2(x_2)$.  The trap confines the condensate well inside $\Omega$,
where $u$ falls below $10^{-6}$ at the boundary, so the TF profile is
compatible with the Dirichlet condition.

The ground state reference solution is computed by the modified $H^1$ flow
on an $800^2$ finite difference grid, run to convergence:
\begin{equation}
  E^* = 3.9053826591, \qquad \lambda^* = 10.1793521806 .
  \label{eq:stirrer-ref}
\end{equation}

\medskip\noindent\emph{How non-separable is the ground state?}  The product
ansatz represents $u=\sqrt{\rho_1}\otimes\sqrt{\rho_2}$, a rank-one matrix on
a tensor grid, so the smallest error any product map can achieve is the
rank-one truncation error of $u^*$.  For~\eqref{eq:stirrer2d} the normalized singular values of
$u^*$ are
$1,\;0.212,\;0.035,\;0.008,\;0.002,\dots$, so
\begin{equation}
  \min_{\rho_1,\rho_2}\frac{\norm{u^*-\sqrt{\rho_1}\otimes\sqrt{\rho_2}}}{\norm{u^*}}
  = 0.211,
  \label{eq:rank1floor}
\end{equation}
which was computed via the singular value decomposition of $u^*$ on the
$800^2$ grid, the floor being
$\bigl(\sum_{k\ge2}\sigma_k^2\big/\sum_{k\ge1}\sigma_k^2\bigr)^{1/2}$.
The same quantity is $0.048$ for the 2D lattice of
Section~\ref{sec:2d_comparison}, whose ground state happens to be close to a
product.  The stirrer is thus about four times more strongly coupled, and
$0.211$ is a floor no product map can cross, however many parameters it is
given.

\medskip\noindent\emph{The full map.}  We replace the product map by a single
boundary-preserving Neural ODE on $\Omega$,
\begin{equation}
  \frac{dw_k}{ds} = \Bigl(1-\frac{w_k^2}{L^2}\Bigr) g_{\theta,k}(w),
  \qquad k=1,2, \qquad w(0)=z, \qquad T_\theta(z)=w(1),
  \label{eq:fullmap}
\end{equation}
where $g_\theta\colon\R^2\to\R^2$, so that each velocity component depends
on both coordinates.  The factor $(1-w_k^2/L^2)$ makes the $k$-th component
vanish on $\{w_k=\pm L\}$, so the Dirichlet condition is preserved exactly as
in Section~\ref{sec:alg:ode}.  The augmented ODE of
Section~\ref{sec:alg:ode} generalizes: writing $A=\partial v/\partial w$
and $J=\partial w/\partial z$ for the Jacobian matrix of the
flow (whose determinant is the scalar $J$ of Section~\ref{sec:alg:ode}),
\begin{equation}
  \dot J = A J, \qquad
  \dot\ell = \operatorname{tr}A, \qquad
  \dot q = J^\top\nabla_w(\operatorname{tr}A),
  \label{eq:aug-full}
\end{equation}
with $J(0)=I$, $\ell=\log|\det J|$ and $q=\nabla_z\ell$, and the score is
recovered as $\nabla_x\log\rho = J^{-\top}(\nabla_z\log\refd - q)$ without
differentiating through the solver.  The state per particle grows from four
components to nine, and the only second derivatives needed are those of the
network with respect to its input, formed analytically.
Both ansatz use the per-axis Beta(5,5) reference density
of~\eqref{eq:beta55}, here on $(-6,6)$, and the solver of
Section~\ref{sec:alg:solver} with $\alpha=0.005$, $N_{\rm ODE}=10$, $200$ CG
iterations and clip $C=10$; the 3D and 4D tests below keep these settings.

\medskip\noindent\emph{Results.}  The product map is stuck $6.9$ to $8.5\%$ above
$E^*$ and does not improve as $M$ grows from $282$ to $1346$, which is the
floor~\eqref{eq:rank1floor} asserting itself, while the full map reaches
$E^*$ within $1$--$2\%$ at every width and with half the parameters
(Table~\ref{tab:stirrer2d}).  Figure~\ref{fig:stirrer2d} shows what the
missing coupling costs: the product map returns a four-lobe pattern where the
ground state has a crescent around the obstacle.  The full map's dense metric
is also the better behaved one, one to two orders of magnitude better
conditioned and at $70$ to $95\%$ of full rank against
$36$ to $45\%$, so the coupling that
a non-separable potential forces on $\hatG$ is not an obstacle.
}

\begin{table}[!ht]
\centering
\caption{2D stirrer~\eqref{eq:stirrer2d}: product map
versus full map at three widths, all with $N=4000$ particles and $600$ PWGF
iterations on one NVIDIA GH200 GPU.  $E_{\rm grid}$ evaluates the reconstructed
density on the $800^2$ grid used for~\eqref{eq:stirrer-ref} and energy err.\
is $|E_{\rm grid}-E^*|/E^*$.  The last two columns describe $\hatG$ at the
returned parameters: $\kappa^+$ is the condition number over the numerically
positive spectrum (eigenvalues above $10^{-14}\lambda_{\max}$)
and $\operatorname{rank}\hatG$ the numerical rank at
relative threshold $10^{-10}$, at most $M$ since $\hatG$ is $M\times M$.}
\label{tab:stirrer2d}
{\small\setlength{\tabcolsep}{5pt}
\begin{tabular}{llrrcrrrr}
\toprule
& & & & energy & & train & & \\
Ansatz & $H$ & $M$ & $E_{\rm grid}$ & err.\ & $\norm{u-u^*}_{L^2}$ & (s)
& $\kappa^+(\hatG)$ & $\operatorname{rank}\hatG$\\
\midrule
Product & 10 & 282  & $4.2128$ & $7.9\%$ & $0.295$ & $107$
        & $6.7\times10^{11}$ & $127$\\
        & 16 & 642  & $4.1735$ & $6.9\%$ & $0.272$ & $116$
        & $1.3\times10^{13}$ & $234$\\
        & 24 & 1346 & $4.2385$ & $8.5\%$ & $0.309$ & $117$
        & $5.1\times10^{12}$ & $525$\\
\addlinespace
Full    & 10 & 162  & $\mathbf{3.9457}$ & $\mathbf{1.0\%}$ & $\mathbf{0.133}$ & $71$
        & $\mathbf{1.7\times10^{10}}$ & $\mathbf{154}$\\
        & 16 & 354  & $3.9513$ & $1.2\%$ & $0.145$ & $78$
        & $2.1\times10^{11}$ & $298$\\
        & 24 & 722  & $3.9702$ & $1.7\%$ & $0.161$ & $75$
        & $6.6\times10^{10}$ & $508$\\
\bottomrule
\end{tabular}
}
\end{table}

\subsubsection{3D stirrer potential}
\label{sec:3d-stirrer}

\begin{figure}[!ht]
  \centering
  \includegraphics[width=\textwidth]{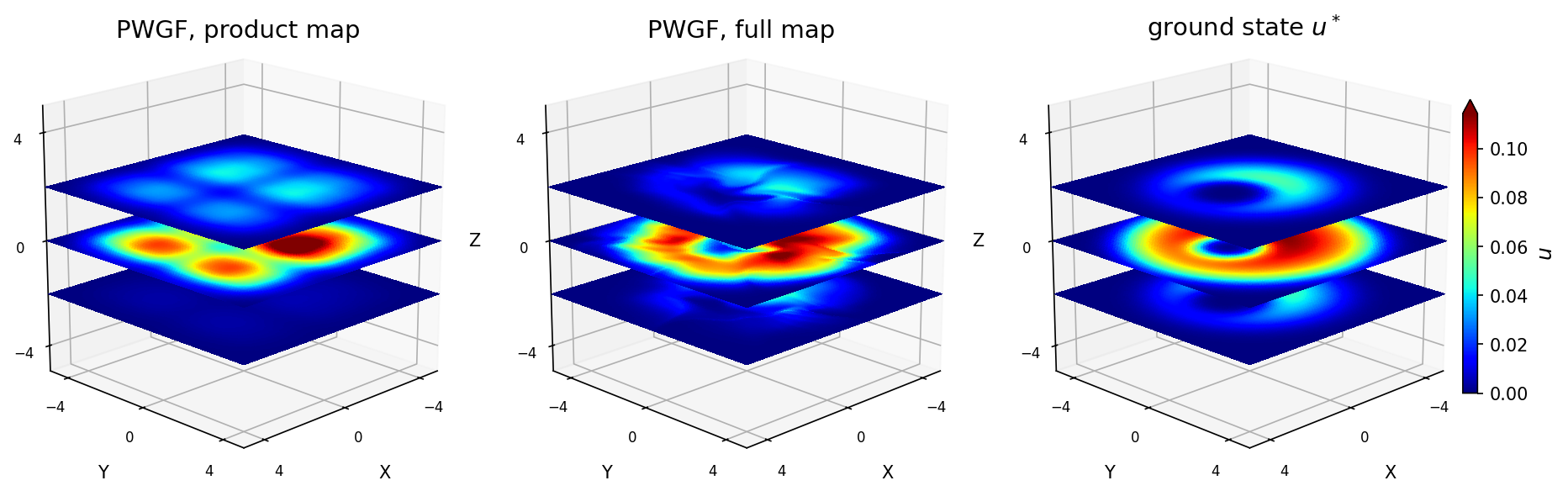}
  \caption{{3D stirrer: slices at $z\approx-2,0,+2$ on a common
    colour scale capped at the peak of the ground state.  \emph{Left:} the
    PWGF product map.  \emph{Middle:} the PWGF full map.  \emph{Right:} the
    ground state $u^*$.  The product map replaces the crescent by four lobes,
    the same failure as in two dimensions (Figure~\ref{fig:stirrer2d}).}}
  \label{fig:stirrer3d-slices}
\end{figure}

\begin{figure}[!ht]
  \centering
  \includegraphics[width=\textwidth]{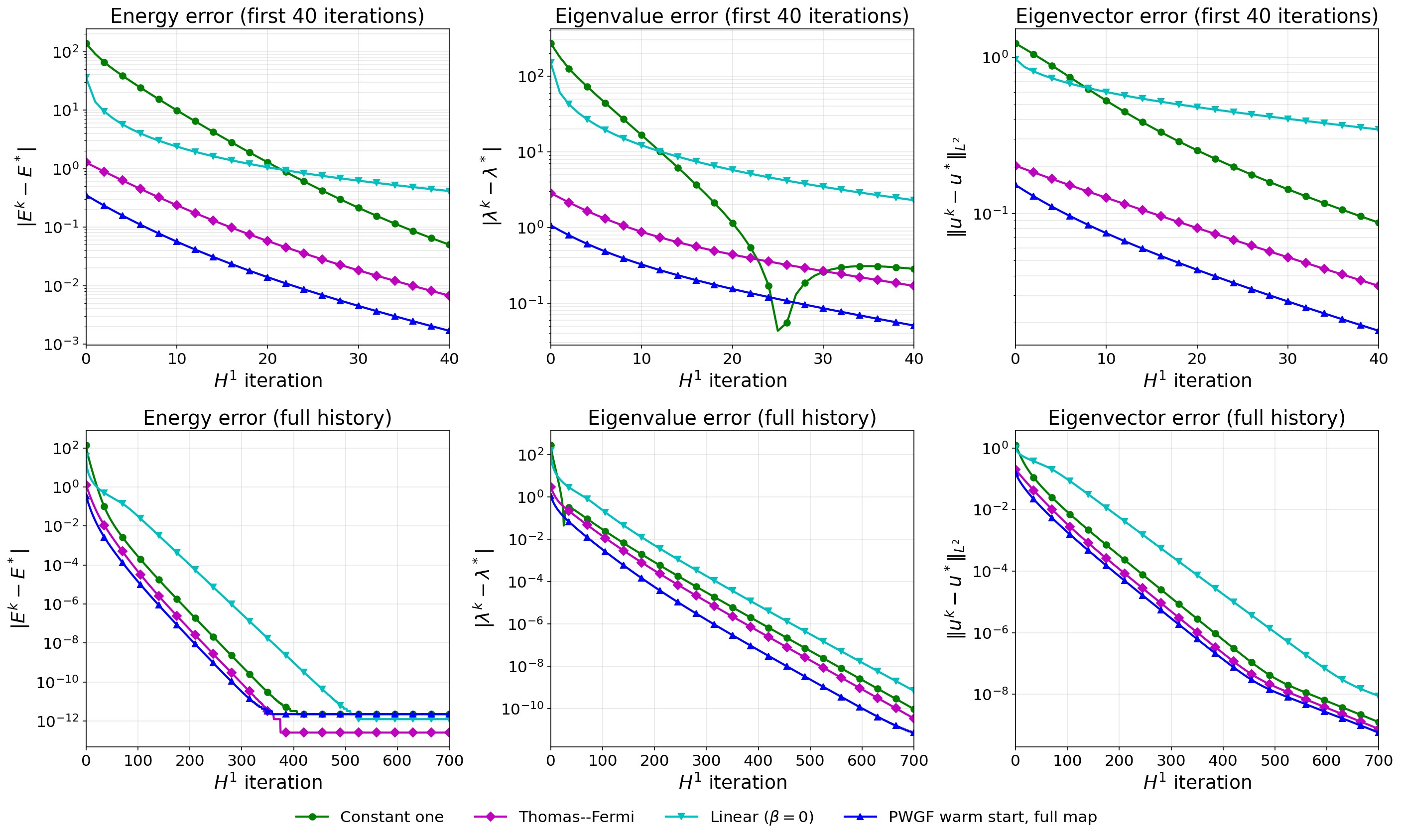}
  \caption{3D stirrer at $999^3$: four initializations under
    the modified $H^1$ flow, double precision, against the ground state
    reference solution~\eqref{eq:stirrer3d-ref}.  \emph{Top row:} the first
    $40$ iterations.  \emph{Bottom row:} all $700$ iterations.}
  \label{fig:stirrer3d-inits}
\end{figure}

{Section~\ref{sec:2d-stirrer} showed that a non-separable potential defeats
the product map and that a full map removes the obstruction.
Section~\ref{sec:3d-warmstart} compared PWGF against standard
initializations at $999^3$.  This subsection combines the two: the same
$999^3$ comparison on a non-separable 3D potential, with both ansatz among
the candidates.  All conventions are those of
Section~\ref{sec:3d-warmstart} unless stated, so that only the potential
changes.

\medskip\noindent\emph{Problem.}  We take the three-dimensional analogue
of~\eqref{eq:stirrer2d}: an anisotropic harmonic trap stirred by a Gaussian
beam directed along $x_3$,
\begin{equation}
  V(\mathbf{x}) = x_1^2 + x_2^2 + \gamma_z^2 x_3^2
  + 2w_0\,e^{-\delta\left((x_1-r_0)^2 + x_2^2\right)},
  \label{eq:stirrer3d}
\end{equation}
with $w_0=30$, $\delta=1$, $r_0=1$, $\gamma_z=2$, $\beta=1600$ on
$\Omega=(-8,8)^3$, the standard model of a laser stirrer in a
condensate~\cite{AntoineDuboscq2015,bao2004}.  The obstacle couples $x_1$
and $x_2$ and is invariant in $x_3$, so the beam bores a channel.

The ground state reference solution is obtained as before, by the modified
$H^1$ flow in double precision to residual stall:
\begin{equation}
  E^* = 9.012310296639, \qquad \lambda^* = 23.358304855681 .
  \label{eq:stirrer3d-ref}
\end{equation}

\medskip\noindent\emph{The cost of separability.}  The smallest eigenvector
error any product map can attain is the rank-one truncation error of $u^*$.
Beyond $d=2$ the singular value decomposition has no exact analogue, so this
is computed by alternating least squares in the mass-weighted $L^2$ norm:
\begin{equation}
  \min_{\rho_1,\rho_2,\rho_3}
  \frac{\norm{u^*-\sqrt{\rho_1}\otimes\sqrt{\rho_2}\otimes\sqrt{\rho_3}}_{L^2}}
       {\norm{u^*}_{L^2}} = 0.278 .
  \label{eq:rank1floor3d}
\end{equation}
The same quantity for the optical lattice of
Section~\ref{sec:3d-warmstart} is $0.186$, computed on the $99^3$ grid in
\ref{sec:3d-quadrature}.  The stirrer thus sits even further from every
product density, and the more consequential difference is qualitative: the
lattice ground state keeps the axis-aligned $3\times3\times3$ blob structure
that a product map can reproduce (Figure~\ref{fig:3d-slices}), whereas no
product density can form the crescent around the obstacle
(Figure~\ref{fig:stirrer3d-slices}).

\medskip\noindent\emph{The two PWGF maps.}  Both ansatz are trained with
$N=4096$ particles for $K=8000$ iterations at width $H=16$
(Table~\ref{tab:stirrer3d-maps}), with the Beta(5,5)$^3$
reference~\eqref{eq:beta55} and the settings of
Section~\ref{sec:2d-stirrer}.  The pattern of
Section~\ref{sec:2d-stirrer} sharpens in three dimensions: the full map
attains the lower energy with $2.5$ times fewer parameters in half the
training time, and its dense metric is over three orders of magnitude better
conditioned and of essentially full numerical rank, against $40\%$ for the
block-diagonal one.

\medskip\noindent\emph{Evaluating the warm start on the $999^3$ grid.}
Both maps are evaluated by backward inversion of the Neural
ODE, with the implicit step of Section~\ref{sec:alg:eval}, so that the
comparison isolates the ansatz.  For the product map the forward tabulation of
Section~\ref{sec:alg:eval} applies as well and serves as an independent check,
agreeing with the inversion to $7\times10^{-5}$ in relative $L^2$.

\medskip\noindent\emph{Results.}  The full-map warm start is the best of the
five, $E^{(0)}=9.365$ against $10.280$ for TF, and keeps the smallest energy
error at every one of the $40$ iterations, ending a factor of $4.0$ below TF
and $30$ below constant-one (Table~\ref{tab:stirrer3d-inits} and
Figures~\ref{fig:stirrer3d-slices} and~\ref{fig:stirrer3d-inits}).  The floor~\eqref{eq:rank1floor3d}
explains why the product map cannot follow: its initial eigenvector error
$0.334$ is within $20\%$ of the floor $0.278$, so the ansatz is essentially
exhausted, while the full map reaches $0.153$, which only coupled
coordinates make possible.  The linear start, which diverges on the optical
lattice, is merely mediocre here, the harmonic trap giving a broad Gaussian
rather than a single-well spike.
}

\begin{table}[!ht]
\centering
\caption{3D stirrer~\eqref{eq:stirrer3d}: the two PWGF
ansatz, trained with $N=4096$ particles for $K=8000$ iterations at $H=16$ on
one NVIDIA GH200 GPU.
$E_{\rm particle}$ is the energy of the returned parameters, estimated by
Monte Carlo on the training particles, and $\kappa^+$ the condition number
over the numerically positive spectrum of $\hatG$
(eigenvalues above $10^{-14}\lambda_{\max}$);
$\operatorname{rank}\hatG/M$ is the numerical rank at relative threshold
$10^{-10}$, as a fraction of $M$.}
\label{tab:stirrer3d-maps}
{\small\setlength{\tabcolsep}{6pt}
\begin{tabular}{lrrrrr}
\toprule
Ansatz & $M$ & $E_{\rm particle}$
& $\kappa^+(\hatG)$ & $\operatorname{rank}\hatG/M$ & train (s)\\
\midrule
Product & $963$ & $9.6665$ & $9.4\times10^{13}$ & $40\%$ & $1848$\\
Full    & $387$ & $\mathbf{8.6786}$ & $\mathbf{1.8\times10^{10}}$ & $\mathbf{99\%}$ & $\mathbf{966}$\\
\bottomrule
\end{tabular}
}
\end{table}

\begin{table}[!ht]
\centering
\caption{3D stirrer~\eqref{eq:stirrer3d} at $999^3$: energy
error after $k$ modified $H^1$ iterations, $1.02$\,s each in single precision
on one NVIDIA GH200 GPU, against the ground state reference
solution~\eqref{eq:stirrer3d-ref} $E^*=9.0123$.  Both PWGF rows are built by backward inversion, so the
comparison isolates the ansatz.  Figure~\ref{fig:stirrer3d-inits}
repeats four of these runs in double precision ($2.49$\,s per iteration);
every tabulated digit is unchanged.}
\label{tab:stirrer3d-inits}
{\small
\begin{tabular}{lrrrr}
\toprule
& \multicolumn{4}{c}{$|E^k-E^*|$, modified $H^1$ flow}\\
Initialization & $k=1$ & $10$ & $20$ & $40$\\
\midrule
Constant one      & $9.1\mathrm{e}{+1}$ & $9.9\mathrm{e}{+0}$ & $1.3\mathrm{e}{+0}$ & $5.0\mathrm{e}{-2}$\\
Linear ($\beta=0$)& $1.4\mathrm{e}{+1}$ & $2.4\mathrm{e}{+0}$ & $1.1\mathrm{e}{+0}$ & $4.1\mathrm{e}{-1}$\\
Thomas--Fermi     & $1.1\mathrm{e}{+0}$ & $2.4\mathrm{e}{-1}$ & $5.8\mathrm{e}{-2}$ & $6.7\mathrm{e}{-3}$\\
PWGF, product map & $8.1\mathrm{e}{-1}$ & $2.3\mathrm{e}{-1}$ & $8.1\mathrm{e}{-2}$ & $1.6\mathrm{e}{-2}$\\
PWGF, full map    & $\mathbf{2.9\mathrm{e}{-1}}$ & $\mathbf{5.6\mathrm{e}{-2}}$ & $\mathbf{1.4\mathrm{e}{-2}}$ & $\mathbf{1.7\mathrm{e}{-3}}$\\
\bottomrule
\end{tabular}
}
\end{table}

\subsection{4D Test}
\label{sec:4d}

\begin{figure}[!ht]
  \centering
  \includegraphics[width=\textwidth]{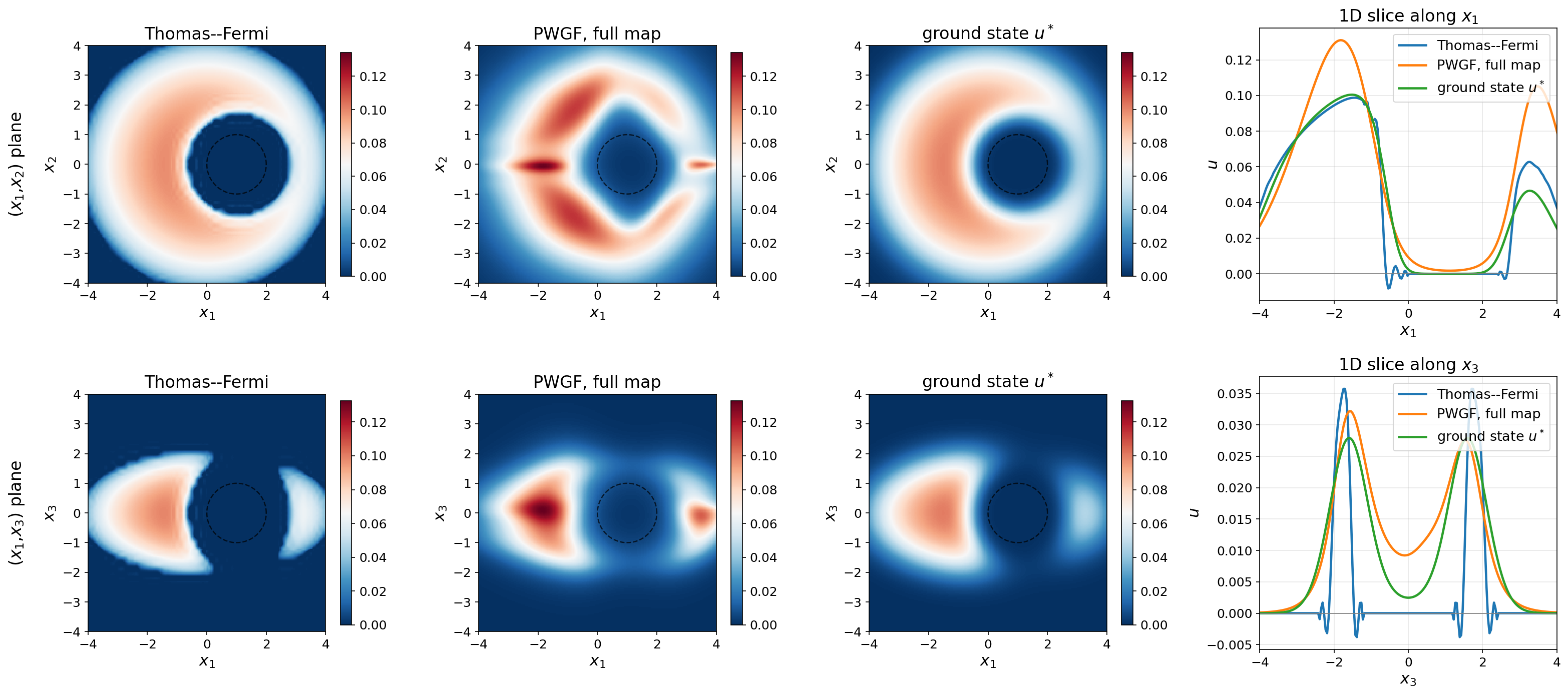}
  \caption{Slices of the TF, the PWGF and the ground state $u^*$ at $99^4$.
    \emph{Top:} $(x_1,x_2)$ plane through the obstacle centre.
    \emph{Bottom:} $(x_1,x_3)$ plane, compressed by the tighter trap
    $\gamma_3=2$.  The dashed circle is the $1/\sqrt\delta$ contour of the
    obstacle and the rightmost column gives 1D slices.   }
  \label{fig:4d-slices}
\end{figure}

\begin{figure}[!ht]
  \centering
  \includegraphics[width=\textwidth]{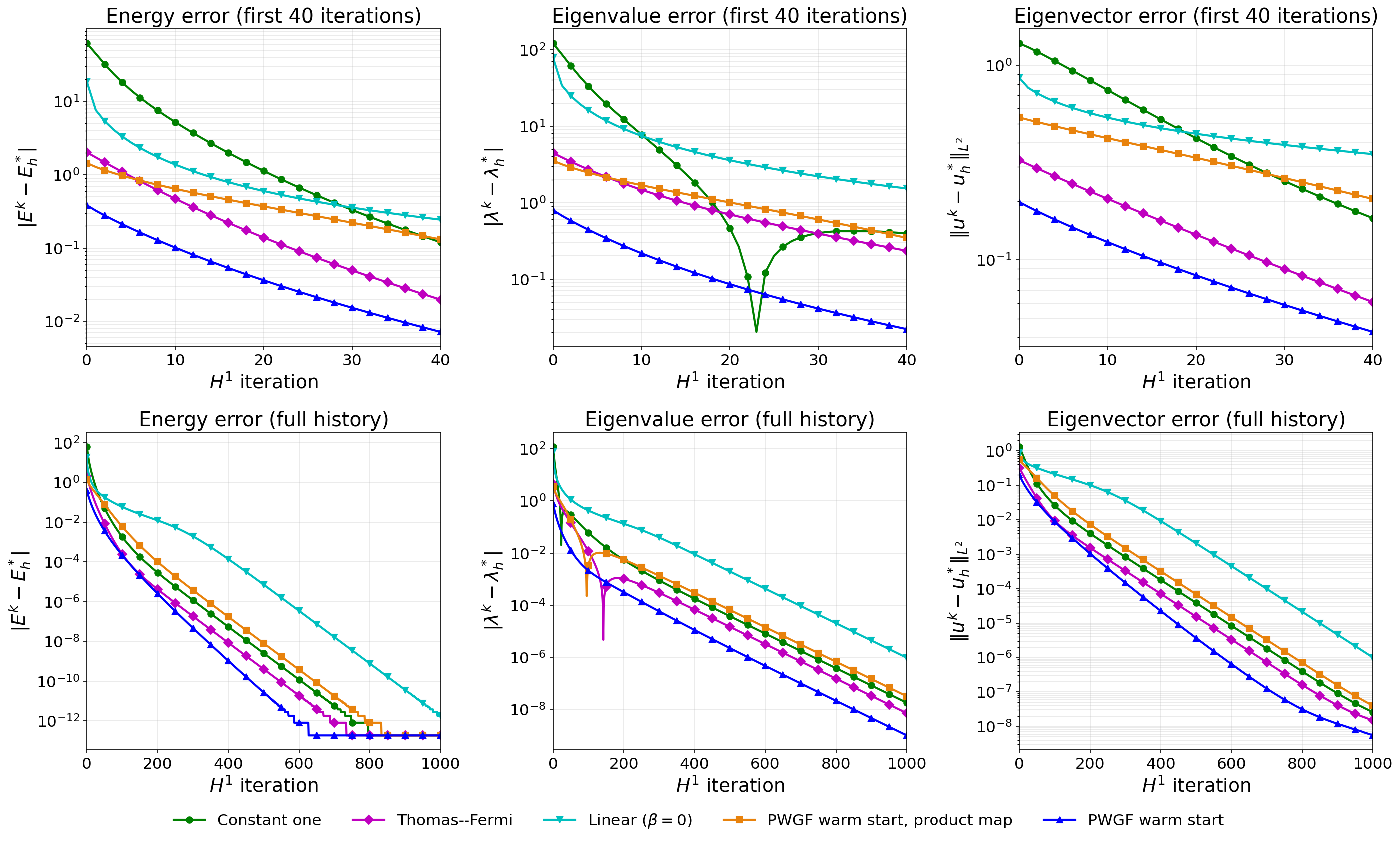}
  \caption{Complete modified $H^1$ history on the $99^4$
    grid from five initializations, double precision, against the
    ground state reference solution~\eqref{eq:ref4d}.  \emph{Top row:} the
    first $40$ iterations,
    with the full-map warm start (blue) lowest throughout.  \emph{Bottom
    row:} $1000$ iterations.}
  \label{fig:4d-inits}
\end{figure}

{The Gross--Pitaevskii equation models a condensate in physical space, so
$d\le3$ exhausts its direct physical content.  Four dimensions are
nonetheless of current interest for superfluids, where vortex cores become
two-dimensional surfaces and the 4D equation has been solved numerically to
study their stability~\cite{mccanna2021}.  High-dimensional Schr\"odinger
problems have also been studied in the neural-network setting through
spectral Barron regularity theory~\cite{chen2023barron}, which gives
approximation rates for two-layer networks that do not
deteriorate as the dimension grows.  It is natural to ask whether a
mesh-free method degrades as the dimension is raised.  We therefore
construct a 4D ground-state problem for the purpose of testing PWGF, and
report what it shows, on cost as well as on accuracy.

\medskip\noindent\emph{Problem.}  We keep the construction of
Section~\ref{sec:3d-stirrer} and change only what four dimensions forces.
The 3D beam becomes a localized bump,
\begin{equation}
  V(\mathbf{x}) = \sum_{k=1}^{4}\gamma_k^2x_k^2
  + 2w_0\,e^{-\delta|\mathbf{x}-\mathbf{r}_0|^2},
  \label{eq:obstacle4d}
\end{equation}
with $\gamma=(1,1,2,2)$, $w_0=120$, $\delta=1$, $\mathbf{r}_0=(1,0,0,0)$ and
$\beta=1600$ on $\Omega=(-6,6)^4$, the same $\beta$ as in the
three-dimensional problems of Sections~\ref{sec:3d}
and~\ref{sec:3d-stirrer}.  With $\norm{u}_{L^2}=1$ the density spreads over a
four-dimensional volume, so $\int u^4$ falls sharply with $d$ and a given
$\beta$ buys much less interaction.  At $\beta=1600$ the interaction is
$19.3\%$ of the energy (kinetic $0.78$, potential $5.93$, interaction
$1.60$), a genuine Gross--Pitaevskii regime.

Discretization is $Q^{10}$ spectral elements with $10$ cells per axis,
$99^4=9.6\times10^{7}$ degrees of freedom.  The separable operator
$-\Delta+V_1$ is inverted by per-axis eigendecompositions and the
non-separable part is carried by preconditioned conjugate gradients, exactly
the 4D solver of~\cite{liu2025gpu}.

\medskip\noindent\emph{Reference.}  No 4D Gross--Pitaevskii ground state is
available in the literature, so we generate one by running the modified $H^1$
flow on the $99^4$ grid in double precision to residual stall,
\begin{equation}
  E^* = 8.318275948409,\qquad \lambda^* = 19.843753363184 .
  \label{eq:ref4d}
\end{equation}

\medskip\noindent\emph{The cost of separability in 4D.}  The rank-one
truncation error of $u^*$, by alternating least squares in the
mass-weighted $L^2$ norm, is
\begin{equation}
  \min_{\rho_1,\dots,\rho_4}
  \frac{\norm{u^*-\bigotimes_k\sqrt{\rho_k}}_{L^2}}{\norm{u^*}_{L^2}} = 0.345 .
  \label{eq:rank1floor4d}
\end{equation}
Against $0.211$ in two dimensions~\eqref{eq:rank1floor} and $0.278$ in
three~\eqref{eq:rank1floor3d}, the floor grows monotonically with
dimension: a product density cannot deplete a localized ball without
depleting $d$ entire hyperplanes.  The growth is a property of the obstacle
strength rather than of the dimension alone.  Carrying the 3D value
$w_0=30$ over gives a floor of only $0.247$, and $w_0$ was raised to $120$
to make the 4D problem the harder one.  Height is the right knob because it
buys coupling almost free of the nonlinearity (the interaction fraction
moves only from $0.203$ to $0.193$), whereas widening the obstacle dilutes
the interaction toward a linear problem (floor $0.133$ at $\delta=0.25$).

\medskip\noindent\emph{The two maps.}  Both ansatz are trained with
$N=4096$ particles for $K=8000$ iterations at $H=16$
(Table~\ref{tab:4d-maps}), with the per-axis Beta(5,5)
reference on $(-6,6)$ and the settings of Section~\ref{sec:2d-stirrer}.  Every effect of two and three dimensions reappears
more strongly: the full map reaches the lower energy with $3.1$ times fewer
parameters in a third of the training time, its dense metric is over three orders
of magnitude better conditioned, and it has full numerical rank where the
block-diagonal metric is at $35\%$.  The redundancy of $d$ independent
copies of one architecture worsens as $d$ grows.

\medskip\noindent\emph{Results.}  The full map begins at $E^{(0)}=8.702$, a
gap of $4.6\%$ against $24.2\%$ for TF, and is lowest in both energy and
$L^2$ at every one of the $40$ iterations
(Figures~\ref{fig:4d-slices} and~\ref{fig:4d-inits}).  The factor of $18$
between the two ansatz at $k=40$ is~\eqref{eq:rank1floor4d} asserting itself:
the product map's initial eigenvector error $0.54$ lies above the floor
$0.345$, the full map's $0.20$ well below.  The ordering persists over the
whole history: the full map reaches an energy error of $10^{-8}$ in $341$
iterations against $395$ for TF and $493$ for the product map
(Table~\ref{tab:h1-iters-4d}).

\medskip\noindent\emph{Cost.}  At this size the offline cost is not repaid.
Training the full map takes $917$\,s, against $0.15$\,s per $H^1$
iteration in double precision on one NVIDIA GH200 GPU.  Treating that
$917$\,s as an offline cost, the warm start saves $12$ iterations against TF,
about $2$\,s of GPU time, at an energy error of $7\times10^{-3}$, and $54$
iterations, about $8$\,s, at $10^{-8}$
(Table~\ref{tab:h1-iters-4d}), so even at its most favourable the saving is
two orders of magnitude below what it cost to obtain.  The trained object is a function
rather than an array, so the balance moves with refinement, as the $999^3$
comparison of Section~\ref{sec:3d-warmstart} shows; locating the 4D
crossover would require far larger grids.

\medskip
In summary, PWGF transfers to four dimensions without modification and is the
best of the five initializations, with the full map dominating the product
map by a wider margin than in three dimensions, as the growth of the
floor~\eqref{eq:rank1floor4d} predicts.  The qualification is cost rather
than accuracy.
}

\begin{table}[!ht]
\centering
\caption{4D obstacle problem~\eqref{eq:obstacle4d}: the two
PWGF ansatz, trained with $N=4096$ particles for $K=8000$ iterations at
$H=16$ on one NVIDIA GH200 GPU.  $E_{\rm particle}$ is the energy of the returned parameters,
estimated by Monte Carlo on the training particles;
$\kappa^+$ and $\operatorname{rank}\hatG/M$ are as in
Table~\ref{tab:stirrer3d-maps}. }
\label{tab:4d-maps}
{\small\setlength{\tabcolsep}{6pt}
\begin{tabular}{lrrrrr}
\toprule
Ansatz & $M$ & $E_{\rm particle}$
& $\kappa^+(\hatG)$ & $\operatorname{rank}\hatG/M$ & train (s)\\
\midrule
Product & $1284$ & $9.4842$ & $1.2\times10^{13}$ & $35\%$ & $2632$\\
Full    & $420$  & $\mathbf{7.9757}$ & $\mathbf{4.2\times10^{9}}$
        & $\mathbf{100\%}$ & $\mathbf{917}$\\
\bottomrule
\end{tabular}
}
\end{table}

\begin{table}[!ht]
\centering
\caption{Modified $H^1$ iterations to reach a given energy
error on the $99^4$ grid.  The PWGF saving is $54$ iterations against TF and
$115$ against constant-one.  The product
map needs $493$ iterations against the full map's $341$ and is slower than
TF despite starting closer in energy, which is the rank-one
floor~\eqref{eq:rank1floor4d} expressed as an iteration count.}
\label{tab:h1-iters-4d}
{\small
\begin{tabular}{lrrrrr}
\toprule
Initialization & $E^{(0)}$ & $10^{-2}$ & $10^{-4}$ & $10^{-6}$ & $10^{-8}$\\
\midrule
Constant one       & $70.15$ & $72$  & $166$ & $306$ & $456$\\
Thomas--Fermi      & $10.33$ & $48$  & $117$ & $246$ & $395$\\
Linear ($\beta=0$) & $26.77$ & $214$ & $412$ & $566$ & $717$\\
PWGF, product map  & $9.76$  & $90$  & $201$ & $343$ & $493$\\
PWGF, full map     & $\mathbf{8.70}$ & $\mathbf{36}$ & $\mathbf{116}$ & $\mathbf{223}$ & $\mathbf{341}$\\
\bottomrule
\end{tabular}
}
\end{table}

\section{Conclusion}
\label{sec:conclusion}

We developed PWGF, a mesh-free method that computes the GPE ground state
by natural gradient descent in the parameter space of a boundary-preserving
Neural ODE, with the Fisher information evaluated analytically via an
augmented ODE.  Numerical tests in dimensions $d=1$ to $4$,
with separable and non-separable potentials, show that the PWGF
energy error is hard to reduce, but that the PWGF output is an effective initial condition for
the $H^1$ Sobolev gradient flow of~\cite{zhang2024sobolev}.
Convergence analysis of the parameterized Wasserstein gradient flow remains
open and all results in this paper are empirical.     

\section*{Acknowledgements}
This research is partially supported by NSF grants DMS-2208518, DMS-2307465, and DMS-2510829.
This work used DeltaAI at the National Center for
Supercomputing Applications (NCSA) through allocation MTH260013 from the
Advanced Cyberinfrastructure Coordination Ecosystem: Services \& Support
(ACCESS) program, which is supported by U.S.\ National Science Foundation
grants \#2138259, \#2138286, \#2138307, \#2137603, and \#2138296.

\section*{Declaration of competing interest}
The authors declare that they have no known competing financial interests
or personal relationships that could have appeared to influence the work
reported in this paper.

\section*{Data availability}
No data was used for the research described in the article.
The Python code for all numerical experiments is available from
the authors upon request.

\section*{CRediT authorship contribution statement}
\textbf{Xiangxiong Zhang:} Software, Formal analysis, Investigation, Visualization,
Writing -- Original Draft, Writing -- Review \& Editing.
\textbf{Haomin Zhou:} Conceptualization, Methodology, Writing -- Review \& Editing.

\section*{Declaration of Generative AI and AI-assisted technologies in the writing process}
During the preparation of this work the authors used Claude (Anthropic)
in order to assist with mathematical discussions, numerical implementation,
and drafting of the manuscript.  After using this tool, the authors reviewed
and edited the content as needed and take full responsibility for the
content of the publication.

\appendix
\setcounter{figure}{0}
\setcounter{table}{0}
\setcounter{equation}{0}
\renewcommand{\thefigure}{A.\arabic{figure}}
\renewcommand{\thetable}{A.\arabic{table}}
\renewcommand{\theequation}{A.\arabic{equation}}
\section{\texorpdfstring{Gauss--Legendre quadrature computations}{Gauss-Legendre quadrature computations}}
\label{app:quadrature}

{This appendix presents a second set of computations in which
the particle estimates are replaced by Gauss--Legendre quadrature,
in one dimension together with a correction of the boundary vanishing
order of the reference density.  These computations identify the origin of
the error floors observed in Section~\ref{sec:numerics} and quantify the
accuracy ceilings of the product-map ansatz.  They also justify the
recommendation of Table~\ref{tab:hyperparams}: Gauss--Legendre quadrature improves
PWGF as a solver, but it does not improve
the warm start for the $H^1$ flow in two and three dimensions.  Every
computation reported in this appendix, including the $99^3$ finite
difference reference of \ref{sec:3d-quadrature}, runs on the CPU of a MacBook Pro
with an Apple M1 processor; no GPU is used.

\subsection{\texorpdfstring{Gauss--Legendre quadrature in 1D}{Gauss-Legendre quadrature in 1D}}
\label{sec:fisher1d}

The ablation study of Section~\ref{sec:ablation} locates the 1D bottleneck
in the particle estimate of $\FQ$ by elimination.  We now identify the
mechanism, verify it quantitatively, and remove it.

\medskip\noindent\emph{The ansatz cannot change the boundary vanishing
order.}  The boundary-preserving field~\eqref{eq:bpode} has fixed points at
$w=\pm1$ (here $L=1$), so the map derivative there is
$T_\theta'(\pm1)=e^{\mp2g_\theta(\pm1)}$, which is finite and nonzero for
every~$\theta$.  By~\eqref{eq:pushforward}, every representable
density $\rho_\theta=(T_\theta)_\sharp\refd$ vanishes at $\partial\Omega$ to
exactly the same order as the reference.  The Beta(2,2) reference
$\refd(z)=\tfrac34(1-z^2)$ vanishes linearly, whereas the true
ground-state density $\rho^*=(u^*)^2$ vanishes quadratically, since
$(u^*)'(\pm1)=\mp\pi/2\neq0$.  For a general potential the quadratic
vanishing follows from the Hopf boundary-point lemma.  Two consequences
follow.  First, $\rho^*$ is not in the ansatz class.  Second, the Fisher
term of every representable density diverges: in reference coordinates the
Fisher integrand behaves, as $|z|\to1$, like
\begin{equation}
  \frac18\,\bigl|\partial_x\log\rho_\theta\bigr|^2\!\circ T_\theta(z)\;
  \refd(z)
  \;=\;\frac{3}{16\,T_\theta'(\pm1)^{2}}\cdot\frac{1}{1-|z|}+O(1),
  \label{eq:fisher-blowup}
\end{equation}
which is not integrable, so $\FQ(\rho_\theta)=+\infty$ for every $\theta$.
The requirement stated below~\eqref{eq:pushforward}, that $\refd$ vanish to
the same order as the target density, is therefore not a technicality.
Violating it is the 1D error floor.  The finite values reported by the
particle estimator exist only because the most extreme of the $N$ samples
sits at distance $\sim N^{-1/2}$ from the boundary, here
$1-\max_i|z_i|\approx5\times10^{-3}$ at $N=3000$.  Truncating the integral
at distance $\delta$ gives
\begin{equation}
  \FQ(\delta)
  \;=\;s(\theta)\,\ln\tfrac1\delta+O(1),
  \qquad
  s(\theta)=\frac{3}{16}\Bigl[T_\theta'(1)^{-2}+T_\theta'(-1)^{-2}\Bigr].
  \label{eq:fisher-slope}
\end{equation}
Both statements are verified numerically: the measured and predicted
slopes of $\FQ(\delta)$ against $\ln(1/\delta)$ are $0.0389$ and
$0.0389$ at the trained parameters, and $0.3806$ and $0.3804$ at
initialization.  During training $s(\theta)$ falls tenfold because the
optimizer inflates $T_\theta'(\pm1)$ from about $1$ to $3.10$, suppressing
the sampled part of the divergence by distorting the map.  This is the
mechanism behind the tail under-sampling bias, and it makes the plateau
unstable: two double-precision re-runs, which change only the sample
realization, move the floor from $0.036$ to $0.065$--$0.066$ and the energy
estimate from $4.187$ to $4.64$, while $3000$-node quadrature at the latter
parameters gives $4.78$ and the exact value is $+\infty$.

\medskip\noindent\emph{Removing the floor.}  The analysis dictates the
cure: match the vanishing order.  We replace the reference by Beta(3,3),
$\refd(z)=\tfrac{15}{16}(1-z^2)^2$, which vanishes quadratically.  Then
every $\rho_\theta$ vanishes quadratically, $\FQ$ is finite on the whole
ansatz class, and $\rho^*$ itself becomes representable, since the monotone
rearrangement from $\refd$ to $\rho^*$ has a derivative bounded away from
$0$ and $\infty$.  Moreover, in
reference coordinates the Fisher integrand equals
$\tfrac12\,u'(T_\theta(z))^2\,T_\theta'(z)$, a smooth function of $z$, so
the $N$ particles can be replaced by $N$ Gauss--Legendre nodes with weights
$\omega_i\refd(z_i)$: deterministic quadrature of spectral accuracy, with
no change to the algorithm, and the $\pm$-symmetric nodes preserve the
sign-symmetry cancellations.  Table~\ref{tab:fisher1d} compares the four
combinations, with every hyperparameter of Table~\ref{tab:hyperparams}
unchanged.

\begin{table}[!ht]
\centering
\caption{{Reference density and integration rule, all other
settings as in Table~\ref{tab:hyperparams}, in double precision except the
first row, which is the single-precision computation of
Figure~\ref{fig:results}.  ``$E$ (train)'' is the minimized objective and
``$E$ (GL)'' re-evaluates the same parameters by $3000$-node quadrature.
$E^*=4.3587$.  For Beta(2,2) the exact energy is $+\infty$ and the GL
column truncates a divergent integral, which is why it exceeds the training
estimate.  Times are wall-clock on the M1 CPU.}}
\label{tab:fisher1d}
\footnotesize
\setlength{\tabcolsep}{3.5pt}
\begin{tabular}{llllll}
\toprule
Reference & Integration & $E$ (train) & $E$ (GL) &
$\|u_\theta-u^*\|_{L^2}$ & time (s) \\
\midrule
Beta(2,2) & MC, FP32 (Fig.~\ref{fig:results}) & $4.187$ & --- & $0.036$ & $26$ \\
Beta(2,2) & MC, FP64              & $4.643$ & $4.784$ & $0.066$ & $45$ \\
Beta(2,2) & GL                   & $4.962$ & $4.962$ & $0.130$ & $45$ \\
Beta(3,3) & MC                   & $4.535$ & $4.401$ & $0.043$ & $46$ \\
Beta(3,3) & GL                   & $4.35814$ & $4.35814$ & $\mathbf{0.0015}$ & $63$ \\
Beta(3,3) & GL, $N_{\rm ODE}{=}40$, $K{=}2000$
                                 & $4.35858$ & $4.35858$ & $\mathbf{0.0009}$ & $1248$ \\
\bottomrule
\end{tabular}
\end{table}

Neither half of the fix works alone.  Quadrature on the divergent reference
makes matters worse ($0.130$ against $0.066$), since it integrates the
divergence more faithfully than Monte Carlo, and the corrected reference
under Monte Carlo improves the floor only to $0.043$.  Together they remove
it: the error drops $24$-fold to $1.5\times10^{-3}$, the energy lands within
$5.6\times10^{-4}$ of $E^*$, the converged $\FQ=1.2334$ agrees with the exact
$\pi^2/8$.  What
remains is the $O(\Delta s)$ bias of the forward-Euler augmented ODE, not
statistics: quartering $\Delta s$ shrinks $E^*-E$ to $1.2\times10^{-4}$
(last row).

The other experiments are affected differently.  The 3D and 4D references
vanish quartically, so $\FQ$ is finite on their ansatz classes, and the 2D
mixture reference~\eqref{eq:gmix_ref} carries the divergence with a
coefficient of order $10^{-5}$ against $\FQ\approx10^{-1}$, invisible at any
realizable sampling depth.  The quadrature replacement is not restricted to
$d=1$: for product maps with product-form references all three energy terms
factor into one-dimensional integrals, as carried out in
\ref{sec:2d-quadrature}.  Genuinely non-separable problems, such as the full
maps of Sections~\ref{sec:2d-stirrer}--\ref{sec:4d}, must retain Monte
Carlo.

With the Beta(3,3) reference and Gauss--Legendre quadrature the energy error
falls to its $5.6\times10^{-4}$ plateau within ten steps instead of stalling
at $0.17$, the wavefunction error converges to $1.5\times10^{-3}$ instead of
$0.036$, and the reconstructed profile is visually indistinguishable
from~$u^*$.

\subsection{\texorpdfstring{Gauss--Legendre quadrature in 2D}{Gauss-Legendre quadrature in 2D}}
\label{sec:2d-quadrature}

\ref{sec:fisher1d} removed the 1D error floor by matching the
boundary vanishing order of the reference and by replacing Monte Carlo
with quadrature.  Both ingredients apply to the 2D problem of
Section~\ref{sec:2d_comparison}, and we test whether they yield the same
improvement.  The envelope of~\eqref{eq:gmix_ref} vanishes linearly, which
reproduces the formal Fisher divergence of \ref{sec:fisher1d}
with the coefficient of order $10^{-5}$ noted there.  The more
consequential ingredient is quadrature.  The product map with the product
reference makes every representable density a product
$\rho_\theta=\rho_1(x_1)\,\rho_2(x_2)$, and $V$ is multiplicatively
separable, so all three energy terms factor into one-dimensional
integrals,
{\small
\begin{equation}
  \FQ=\sum_{k=1}^{2}\tfrac18\int\bigl[(\log\rho_k)'\bigr]^2\rho_k dx_k,
  \,
  \FV=\prod_{k=1}^{2}\int\sin^2\!\bigl(\tfrac{\pi x_k}{4}\bigr)\rho_k dx_k,
  \,
  \FR=\tfrac{\beta}{4}\prod_{k=1}^{2}\int\rho_k^2 dx_k,
  \label{eq:2d-factorized}
\end{equation}}
and the entire 2D energy is computable by two one-dimensional
Gauss--Legendre rules with $2000$ nodes per axis, with no Monte Carlo
sampling.

To interpret the results, we first compute the accuracy ceiling of the
ansatz.  Every representable $u_\theta=\sqrt{\rho_1}\otimes\sqrt{\rho_2}$
is a rank-one tensor, while the multiplicative coupling makes $u^*$ not a
product (see also the rank-one analysis of Section~\ref{sec:2d-stirrer}).
The singular value decomposition of the ground state $u^*_h$ therefore
gives a floor that no product map can cross:
\[
  \min_{\mathrm{rank\mbox{-}one}}\|u^*_h-u_1\otimes u_2\|_{L^2}=0.048,
  \qquad \sigma_2/\sigma_1=4.7\times10^{-2},
\]
and the normalized best rank-one truncation has grid energy $0.2227$.
The warm start of Section~\ref{sec:h1flow}, with $\|u-u^*_h\|=0.158$ and
$E_{\rm FD}=0.2327$, sits three times above this floor.

Table~\ref{tab:2d-quadrature} reports three double-precision quadrature
runs against the Monte Carlo computation of Section~\ref{sec:h1flow}: the
two envelope orders at the budget of Table~\ref{tab:hyperparams} and the
order-matched envelope run to convergence (at $K=2000$ the last $300$
steps move $E$ by $10^{-5}$).
Figure~\ref{fig:2d-quadrature-history} shows the corresponding histories.

\begin{table}[!ht]
\centering
\caption{{The fix of \ref{sec:fisher1d} applied to the 2D
problem, all other settings as in Table~\ref{tab:hyperparams}, double
precision, $2000$ Gauss--Legendre nodes per axis.  ``Training time'' is the
PWGF wall clock, ``$E$ (FD)'' and $\|u_\theta-u^*_h\|$ evaluate the
reconstructed warm start on the $200^2$
grid, and the $H^1$ column counts iterations to residual $<10^{-12}$
($a=0.15$, $dt=1$), with constant-one taking $81$.  ``Envelope'' is the
boundary factor of~\eqref{eq:gmix_ref}, linear as published or quadratic to
match the vanishing order of $\rho^*$.}}
\label{tab:2d-quadrature}
\small
\setlength{\tabcolsep}{5pt}
\begin{tabular}{llllll}
\toprule
Envelope & Integration & training time (s) & $E$ (FD) &
$\|u_\theta-u^*_h\|$ & $H^1$ iters \\
\midrule
linear    & MC, $K{=}400$  & $78$   & $0.2327$ & $0.158$ & $99$ \\
quadratic & GL, $K{=}400$  & $152$  & $0.2333$ & $0.145$ & $101$ \\
linear    & GL, $K{=}400$  & $192$  & $0.2344$ & $0.164$ & $106$ \\
quadratic & GL, $K{=}2000$ & $1126$ & $0.2302$ & $0.142$ & $108$ \\
\bottomrule
\end{tabular}
\end{table}

\begin{figure}[!ht]
  \centering
  \includegraphics[width=\textwidth]{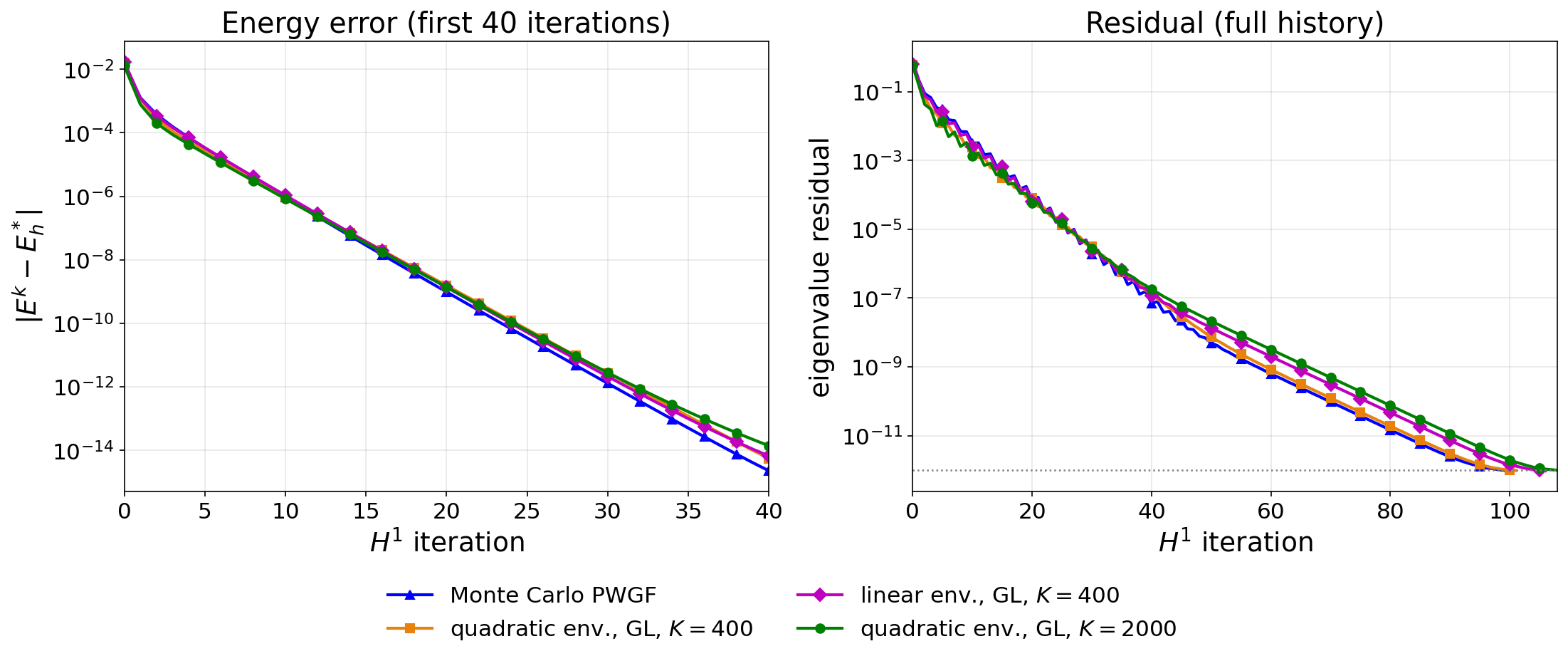}
  \caption{Modified $H^1$ history on the $200^2$ grid from
    the four warm starts of Table~\ref{tab:2d-quadrature}, against
    $E^*_h=0.21706$.  \emph{Left:} the energy error over the first $40$
    iterations, where the four starts lie within a factor of $1.3$ at $k=0$
    and are indistinguishable by $k=10$.  \emph{Right:} the eigenvalue
    residual, which the stopping rule watches (dotted line at $10^{-12}$);
    the curves are parallel and their small offsets are what produce the
    counts $99$--$108$.  The iteration count therefore measures the tail of
    the residual rather than the quality of the start.}
  \label{fig:2d-quadrature-history}
\end{figure}

The estimator bias is removed: the training objective now
agrees with the grid energy to $2\times10^{-4}$, whereas the Monte Carlo
objective $0.2162$ sits $0.017$ below its grid energy $0.2327$.  The
accuracy nevertheless barely moves: $0.158\to0.142$ in $L^2$ and
$0.2327\to0.2302$ in energy, still far from the rank-one floor.  With the
Gauss--Legendre estimator and the budget converged, what remains is the product class
itself.  The 1D floor was created by the
estimator and the fix removed it; the 2D floor is created by the
representation, which no estimator can remove.  The warm start is unaffected
either way, its $H^1$ iteration count even anti-correlating
with its energy (the lowest-energy start takes $108$ iterations, the Monte
Carlo one $99$), so we retain the Monte Carlo setup of
Table~\ref{tab:hyperparams}.  The three-dimensional analogue, with both
ceilings computed exactly, is in \ref{sec:3d-quadrature}.

\subsection{\texorpdfstring{Gauss--Legendre quadrature in 3D}{Gauss-Legendre quadrature in 3D}}
\label{sec:3d-quadrature}

The quadrature construction of \ref{sec:2d-quadrature} applies to
the 3D problem of Section~\ref{sec:3d} in a simpler form, since $V$ is
additively separable.  With the product map and the Beta(5,5)$^3$ product
reference,
{\small
\begin{equation}
  \FQ=\sum_{k=1}^{3}\tfrac18\int\bigl[(\log\rho_k)'\bigr]^2\rho_k\,dx_k,
  \;
  \FV=\tfrac12\sum_{k=1}^{3}\int v(x_k)\,\rho_k\,dx_k,
  \;
  \FR=\tfrac{\beta}{4}\prod_{k=1}^{3}\int\rho_k^2\,dx_k,
  \label{eq:3d-factorized}
\end{equation}}
with $v(x)=x^2+100\sin^2(\pi x/4)$, and three one-dimensional
Gauss--Legendre rules with $2000$ nodes per axis replace the $N=6000$
Monte Carlo particles.  No divergence needs fixing here: the quartic
envelope gives a finite Fisher integral on the whole ansatz class, and the
true density is exponentially small at the walls.

In three dimensions both ceilings of the product ansatz can be computed
exactly.  By symmetry the optimal product state is $u\otimes u\otimes u$,
where $u$ minimizes $3\bigl[\tfrac12\int(u')^2+\tfrac12\int v\,u^2\bigr]
+\tfrac{\beta}{4}\bigl(\int u^4\bigr)^{3}$ over $\|u\|_{L^2}=1$, i.e.,
solves the one-dimensional eigenproblem
\begin{equation}
  -u'' + v\,u + \beta c^2 u^3=\lambda u,
  \qquad c=\int u^4\,dx,
  \label{eq:prod-selfconsistent}
\end{equation}
with the self-consistent constant $c$.
Solving~\eqref{eq:prod-selfconsistent} on the $99^3$ grid gives
\[
  E_{\rm prod}^*=34.504
  \;\;(2.2\%\ \text{above}\ E^*_h=33.769),
  \qquad
  \|u_{\rm prod}-u^*_h\|_{L^2}=0.213,
\]
while the best rank-one approximation of $u^*_h$, computed by alternating
least squares, has error $0.186$ and energy $34.97$.  Every product-map
method, mesh-free or not, is bounded by these numbers.

Table~\ref{tab:3d-quadrature} compares the Monte Carlo PWGF map of
Table~\ref{tab:h1-iters-n99} with quadrature-trained maps and with the
product optimum, and Figure~\ref{fig:3d-quadrature-history} shows the
corresponding energy histories.  The step-size guards of Table~\ref{tab:hyperparams}, tuned for noisy Monte
Carlo gradients, do not stabilize the Gauss--Legendre objective, and a tighter trust
region is needed before it descends cleanly.  The Monte
Carlo bias also acts as annealing: at the shared budget $K=400$, Monte
Carlo training had reached a grid energy of $39.19$ while the
Gauss--Legendre map still evaluated to $52$ on the grid ($47$ in its own
objective), because the biased objective underweights the
Fisher barriers and speeds up early mass transport even though its minimum
is displaced.

\begin{table}[!ht]
\centering
\caption{{Gauss--Legendre quadrature in 3D: warm-start quality on the
$99^3$ grid and modified-$H^1$ iterations to $|E^k-E^*_h|<$ tol, with
constant-one taking $141$/$203$ (Table~\ref{tab:h1-iters-n99}).
``Training time'' is the PWGF wall clock and ``$E$ (FD)'' the grid
evaluation of the reconstructed state.
The first row is the Monte Carlo map of Table~\ref{tab:h1-iters-n99}.}}
\label{tab:3d-quadrature}
\footnotesize
\setlength{\tabcolsep}{3pt}
\begin{tabular}{llllll}
\toprule
Warm start & training time (s) & $E$ (FD) & $\|u_\theta-u^*_h\|$ &
$10^{-6}$ & $10^{-8}$ \\
\midrule
Monte Carlo PWGF & $\approx240$ & $39.19$ & $0.458$
  & $\mathbf{99}$ & $\mathbf{144}$ \\
quadrature, trust region, $K{=}4000$ & $2661$ & $35.74$ & $0.266$
  & $132$ & $193$ \\
MC map $+$ 3000 quad.\ steps & $\approx2410$ & $36.02$ & $0.267$
  & $121$ & $182$ \\
\bottomrule
\end{tabular}
\end{table}

\begin{figure}[!ht]
  \centering
  \includegraphics[width=\textwidth]{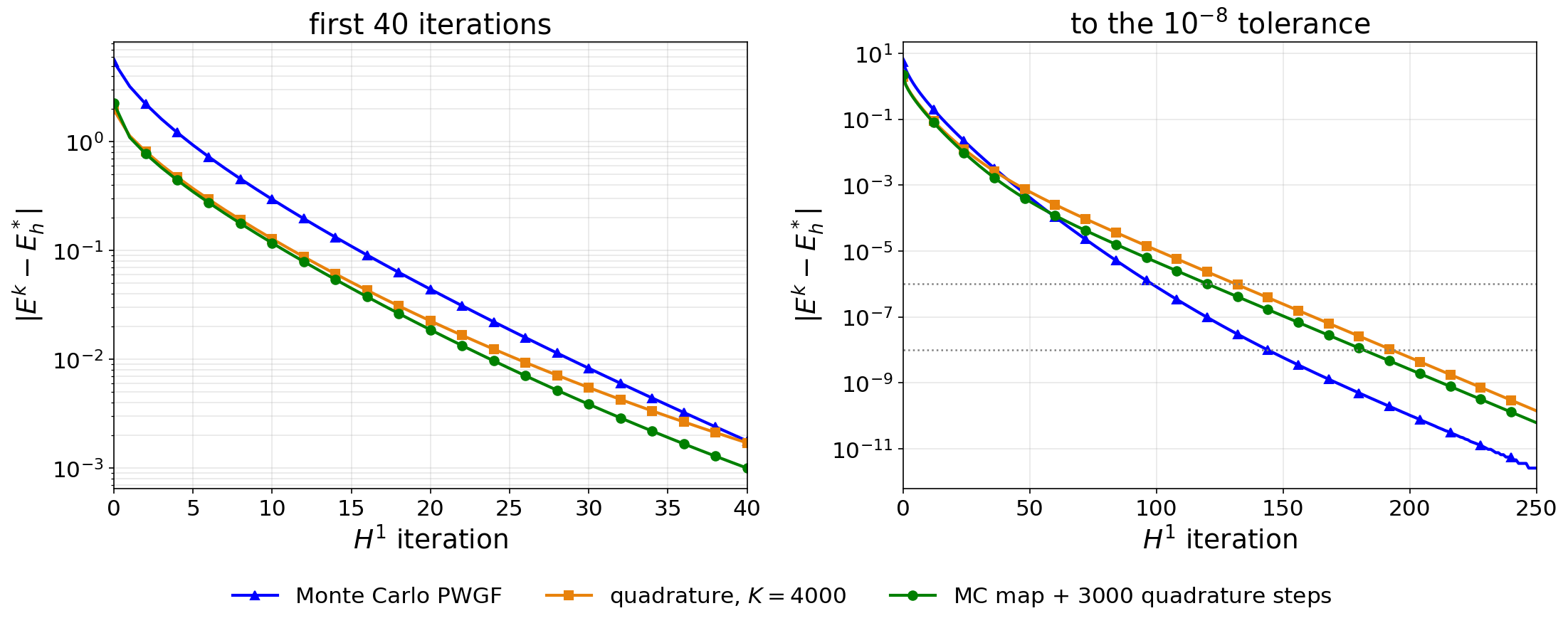}
  \caption{Modified $H^1$ energy history on the $99^3$ grid
    from the three warm starts of Table~\ref{tab:3d-quadrature}, against
    $E^*_h=33.769$.  \emph{Left:} the first $40$ iterations.  \emph{Right:}
    out to the $10^{-8}$ tolerance, with the two tolerances of
    Table~\ref{tab:3d-quadrature} dotted.  The Monte Carlo start begins
    highest, at $5.42$ against $1.97$ and $2.25$, but descends faster,
    crosses the other two at iterations $41$ and $58$ and stays below them,
    which is why it reaches every tolerance first.}
  \label{fig:3d-quadrature-history}
\end{figure}

Gauss--Legendre quadrature buys a better solution rather than a better warm start:
PWGF with Monte Carlo quadrature is already a cheap and effective warm
start, and none of the more accurate maps beats it.

}

\end{document}